\documentclass[11pt]{article}
\usepackage{amsfonts}
\usepackage{amsmath,amsthm,amscd,amssymb,mathrsfs,setspace, textcomp}

\usepackage{latexsym,epsf,epsfig,float}
\usepackage{color}
\usepackage[hmargin=2.25cm,vmargin=2.75cm]{geometry}

\usepackage{hyperref}

\usepackage{array,multirow}
\usepackage{multicol}
\usepackage{cite}
\usepackage{empheq}

\usepackage{enumitem,mathtools}
\usepackage{multimedia}
\usepackage{tikz}
 
\usepackage{pgfplots}
\usepackage{tikz}
\usetikzlibrary{patterns,matrix,arrows,decorations,calc,decorations.pathreplacing,backgrounds,shapes}
\usetikzlibrary{decorations.pathmorphing,decorations.text}
\usetikzlibrary{scopes,decorations.markings,positioning,shapes.arrows} 
\usepackage{graphicx}  
\usepackage{amsmath,amssymb}
\usepackage{bm}
\usepackage{array}
\usepackage{geometry}  
\usepackage{booktabs} 
\usepackage{subfig}
\usepackage{tabularx}
\usepackage{caption}
\usepackage[english]{babel}

\newcommand{\ds}{\displaystyle}

\newcommand{\cA}{{\mathcal{A}}}

\theoremstyle{plain}
\newtheorem{theorem}{Theorem}[section]

\newtheorem{definition}{Definition}

\theoremstyle{remark}
\newtheorem{remark}{Remark}[section]
\numberwithin{equation}{section}
\numberwithin{theorem}{section}
\numberwithin{remark}{section}
\numberwithin{assumption}{section}
\numberwithin{condition}{section}

\title{Large Deflections of Inextensible Cantilevers: \\ Modeling, Theory, and Simulation}

 \author{\normalsize \begin{tabular}[t]{c@{\extracolsep{.6em}}c@{\extracolsep{.6em}}c@{\extracolsep{.6em}}c}
      Maria Deliyianni\footnote{University of Maryland, Baltimore County, MD}&  Varun Gudibanda \footnote{Carnegie Mellon University}& Jason Howell \footnote{Carnegie Mellon University}&  Justin T. Webster\footnote{University of Maryland, Baltimore County, MD}\\
\it mdeliy1@umbc.edu  ~~& \it  vgudiban@andrew.cmu.edu ~~ &  \it  ~~howell4@cmu.edu  & \it websterj@umbc.edu
\end{tabular}}

\begin{document}
\maketitle

\begin{abstract} {\noindent A recent large deflection cantilever model is considered. The principal nonlinear effects come through the beam's {\em inextensibility}---local arc length preservation---rather than traditional extensible effects attributed to fully restricted boundary conditions. Enforcing inextensibility leads to: {\em nonlinear stiffness} terms, which appear as quasilinear and semilinear effects, as well as {\em nonlinear inertia} effects, appearing as nonlocal terms that make the beam implicit in the acceleration. 

In this paper we discuss the derivation of the equations of motion via Hamilton's principle with a Lagrange multiplier to enforce the {\em effective inextensibility constraint}. We then provide the functional framework for weak and strong solutions before presenting novel results on the existence and uniqueness of strong solutions. A distinguishing feature is that the two types of nonlinear terms prevent independent challenges: the quasilinear nature of the stiffness forces higher topologies for solutions, while the nonlocal inertia requires the consideration of Kelvin-Voigt type damping to close estimates. Finally, a modal approach is used to produce mathematically-oriented numerical simulations that provide insight to the features and limitations of the inextensible model. 
  \\[.3cm]
\noindent {\bf Key terms}: nonlinear beam, cantilever, inextensibility, large deflections, quasilinearity
 \\[.3cm]
\noindent {\bf MSC 2010}: 74B20, 35L77, 37L05, 35B65, 70J10}
\end{abstract}

\maketitle

\section{Introduction}

This paper considers a recent partial differential equation (PDE) model for the large deflections of a clamped-free, elastic beam (a {\em cantilever}). Motivated by aeroelastic applications described below, we consider, physically, a thin, narrow plate with an aspect ratio such that the large deflections predominantly exhibit 1-D features---though future work will address fully 2-D plate models. The cantilever model of interest is distinguished by its derivation from an inextensibility constraint: the enforcement of arc-length preservation. This inextensible cantilever model was recently derived in \cite{inext1}, though inextensibility has been treated in a similar fashion for the past 30 years \cite{semler2,paidoussis,inext2}. Enforcing inextensibility in the beam leads to both nonlinear stiffness effects, as well as nonlinear inertial effects. The former yields quasilinear and semilinear terms in the equation of motion, and the latter contributes nonlocal terms that prevent the equation from being written as a traditional second-order-in-time evolution. Grappling with these (independent) nonlinear effects is at the heart of the mathematical and numerical challenge for large deflection cantilever dynamics.

In the engineering literature, as well as the PDE and control literature, beam theory is well studied. Mathematically, the linear theory of Euler-Bernoulli, Rayleigh, shear, and Timoshenko beams---across all boundary configurations---has been established for some time (see, for instance, the nice survey \cite{beams}). Nonlinear beam models, such as Kirchhoff or Krieger-Woinowsky beams, have been considered, typically based on the property of extensibility---see \cite{ball,dickey,HolMar78,wonkrieg} for some older references, as well as the more recent \cite{HTW,HHWW,beam4}. Extensible beams are characterized by a nonlinear restoring force that accounts for the effects of stretching on bending; these are often cubic-type semilinear models. A clear extensible modeling discussion is given in \cite{lagleug}, where a nonlinear beam system (accounting for in-axis and out of axis dynamics) is studied from a semigroup and boundary control point of view, permitting the possibility of the cantilevered configuration. (This model is the beam equivalent of the so called {\em full von Karman} plate model---see \cite{koch}.) Based on the model in \cite{lagleug}, the paper \cite{HTW} studied the well-posedness and long-time behavior of a reduced, scalar version in the cantilevered configuration with a nonconservative loading. Further numerical work appearing in \cite{HHWW} addressed unstable extensible beams across all physical configurations. 

To the knowledge of the present authors, no mathematical theory of solutions, akin to the above references, for the nonlinear inextensible beam has been attempted, and the body of simulations performed using this model have appeared strictly in the engineering literature. Therefore we: \begin{itemize} \setlength\itemsep{.01cm} \item recall the derivation of the equations of motion, \item prescribe a functional setup for the dynamics, \item present the first well-posedness results for strong solutions, \item give a discussion of the associated energy estimates and construction techniques,
\item and provide mathematically-oriented numerical investigations of the dynamics. \end{itemize}

Cantilever models are often utilized in situations where a dynamic driver appears through one of the boundary conditions, or via some distributed forcing function. The primary motivation in this paper comes from aeroelasticity \cite{dowell}, where non-conservative distributed forcing represents aerodynamical pressure differentials across the beam.  Follower forces---maintained tangential forces at the free end---have also recently studied numerically \cite{follower,langre2}.  These non-conservative terms can lead to structural bifurcation and associated large deflections through limit cycles oscillations (LCOs) or even chaos \cite{HolMar78,HHWW}. More specifically, cantilevers in an axial\footnote{the unperturbed flow runs {\em along} the principal axis, as opposed to the more common {\em normal} configuration, where the flow is orthogonal to the beam's span \cite{amjad,balamjad}} airflow can experience an instability known as {\em flutter}, even at low flow velocities. Beyond critical flow parameters, the system enters an LCO of persistent, flapping motions that can be on the order of the beam's length \cite{dowell4,inext2}. It has been shown, for instance, that such dynamics can generate power from which energy can be {\em harvested} \cite{DOWELL,energyharvesting}. To effectively and efficiently do this {one must understand the qualitative properties of the LCO} \cite{DOWELL,experimental}, and hence a proper PDE analysis must treat a nonlinear, large-deflection cantilever model.

In order to provide some visual context, Figure 1 above shows temporal snapshots of a cantilever LCO in actual wind-tunnel experiments. Figure 2 shows simulated snapshots of the first and second Euler-Bernoulli cantilever modes (in vacuo eigenfunctions).

\begin{figure}[htp]
\begin{center}
\includegraphics[width=2in]{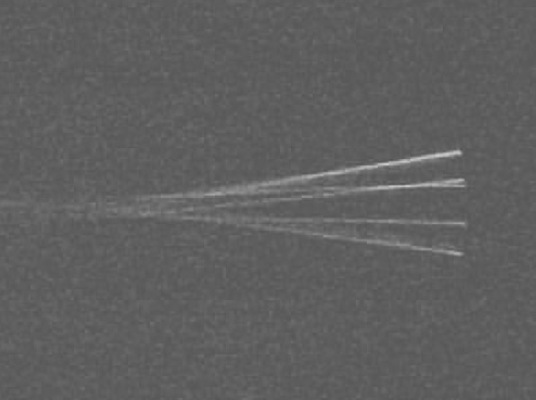}
\hspace*{0.55in}
\includegraphics[width=2in]{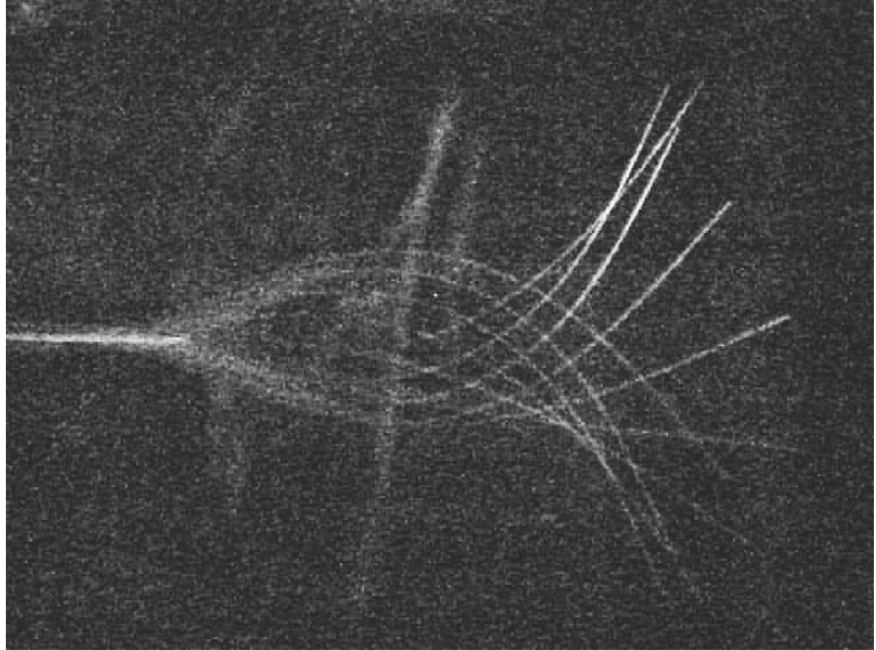}
\vskip-.15cm
\caption{Temporal snapshots of post onset high-velocity LCO (left) and low-velocity LCO (right) for a cantilever. Captured from wind-tunnel simulations \cite{dowell4,inext2}.}
\end{center} 
\end{figure}

\begin{figure}[htp]
\begin{center}
\includegraphics[width=2in]{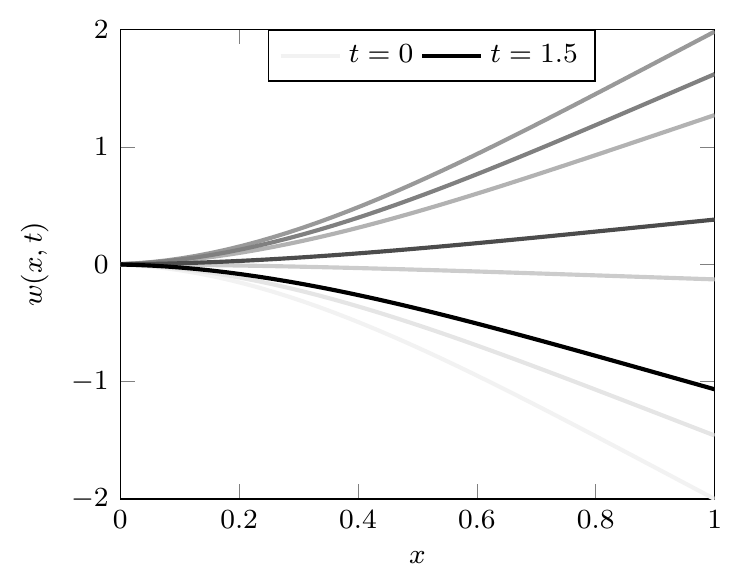}
\hspace*{0.2in}
\includegraphics[width=2in]{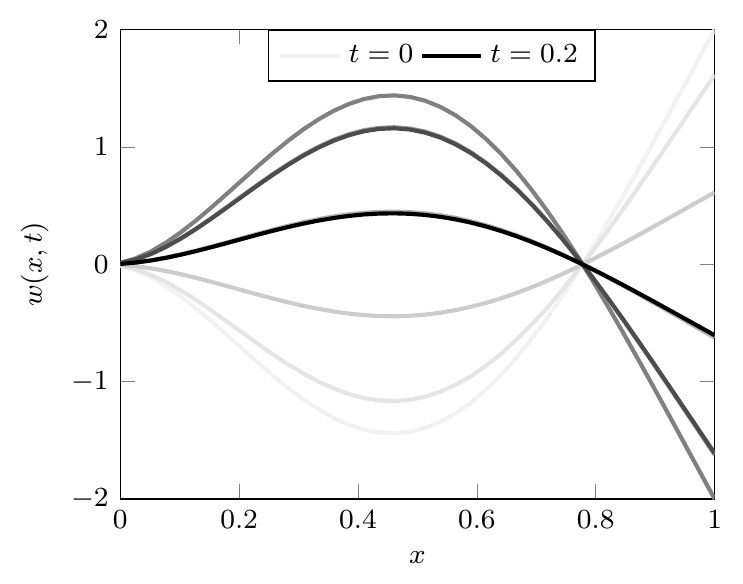}
\vskip-.2cm
\caption{In vacuo \emph{linear} dynamics; temporal snapshots of the the first two Euler-Bernoulli cantilever modes (left) and (right).}
\end{center} 
\end{figure}

\subsection{Applications and Background}

The principal fluid-structure phenomenon associated with nonlinear cantilever deflections is that of {\em aeroelastic flutter:} structural self-destabilization brought about by a surrounding flow. Flutter occurs in many scenarios: elastic structures in wind; aircraft components \cite{dowell}; pipes conveying fluid \cite{paidoussis}; and in human respiration \cite{huang}. From a design point of view, it cannot be overlooked due to large amplitude structural response.   Until about 15 years ago, interest in cantilevers in axial flow
had been minimal  \cite{huang}\footnote{Here we make the distinction with flapping flags---not typically modeled via fourth order equations with stiffness---which have been a topic of interest for hundreds of years; see \cite{flag}.}. On the other hand, interest in airfoil and panel flutter has been {\em immense} for 75 years (see the monograph \cite{dowell}).  In the most prominent cases, flutter is undesirable, with a design goal to prevent it. 

Recently, the axial flow flutter of a beam or plate has been a topic of great interest in the engineering literature \cite{paidoussis, dowell4, TP0,langre2}, as well as, very recently, the mathematical literature \cite{amjad,HTW}. This interest is predominantly due to piezoelectric energy harvesting applications \cite{shubov,piezomass}. The general idea for large displacement harvesters, realized in recent experiments \cite{DOWELL,experimental,tang3}, is to capture mechanical energy in LCOs via piezoelectric laminates or patches (for which oscillating strains induce current \cite{energyharvesting,ozer,scott1}).

As with all {\em flutter problems}, the {\em onset} of instability can be studied from the point of view of a linear structural theory \cite{harvest*1,huang}---typically as an eigenvalue problem (see \cite{shubov,vedeneev,HHWW} for recent discussions). However, if one wishes to study dynamics in the post-flutter regime, the analysis will require {\em some} nonlinear restoring force that will keep solutions {\em bounded in time} \cite{survey2,HTW}.  From \cite{harvest*1}: ``To assess the amount of electrical power that can efficiently be extracted, nonlinear effects are important to provide the saturation amplitude of the self-sustained oscillations." The choice of nonlinear restoring force(s) dictate(s) the qualitative LCO properties. Additionally, power generation considerations require predicting total LCO energy to determine extractable electrical power \cite{piezomass,energyharvesting,tang3,DOWELL}. 
		  Thus, an understanding of nonlinear cantilever deflections, determined from intrinsic parameters, is of critical importance for energy harvesting applications. For this, we must have a robust cantilever model, accommodating the placement of sensors/actuators/piezo devices \cite{scott1}, which translates to spatially localized inertia, stiffness, and damping effects \cite{inext1,scott1}.
		   See also the papers \cite{paidoussis,paidoussis3} for further discussions of the need for, and effects of, implementing the inextensibility constraint in the context of tubes conveying fluid.

 With the above applications in mind, {\em the central challenge is to capture, analyze, and predict cantilever large deflections}; this translates to a viable theory of PDE solutions for the inextensible clamped-free beam, as well as robust and efficient associated computational methods. 
 
 \subsection{Equations of Motion}
 
 We relegate a modeling discussion to the section that follows, but to conclude the introduction, let us state the equations of interest. Consider the (non-dimensionalized) quantities of interest:
 \begin{itemize}  \setlength\itemsep{.01cm}
\item $D>0$,~$L>0$: beam stiffness and length, resp.;
\item $k_2 \ge 0$:  $k_2$ is Kelvin-Voigt type \cite{che-tri:89:PJM,beamdamping} damping\footnote{We choose $k_2$ to denote Kelvin-Voigt type damping, since, in general standard beam damping could be written $\mathbb D w_t=[k_0-k_1\partial_x^2+k_2\partial_x^4]w_t$, indicating weak (frictional), square root-like, and strong damping, respectively. See Section \ref{damping}.}
\item $p(x,t)$ is the distributed loading across the beam span.
\end{itemize}
Now, let $u: [0,L] \times [0,T] \to \mathbb R$ and $w: [0,L] \times [0,T] \to \mathbb R$ correspond (respectively) to the in-axis (longitudinal) and out-of-axis (transverse) Lagrangian deflections.
Then the dynamic equations of motion for the inextensible cantilever are:
\begin{equation}\label{dowellnon}
\begin{cases} \displaystyle w_{tt}+D\partial_x^4w+ k_2 \partial_x^4 w_t+\mathbf A (w) =p(x,t)& \text{ in } (0,L) \times (0,T)
\\ w(t=0)=w_0(x), w_t(t=0)=w_1(x) \\ w(x=0)=w_x(x=0)=0;~w_{xx}(x=L)=w_{xxx}(x=L)=0. 
\end{cases}
\end{equation} 
The nonlinear, nonlocal operator $\mathbf A$ given by
\begin{align}\label{dowellnon2}
\mathbf A (w) =&- D\partial_x\big[(w_{xx})^2w_x \big]+  D\partial_{xx}\big[w_{xx}\big(w_x^2\big)\big] +  \partial_x\left[w_x\int_x^L u_{tt}(\xi)d\xi \right] \\
u(x)=&~-\frac{1}{2}\int_0^x \left [w_x(\xi) \right ]^2d\xi. \label{dowellnon3}
\end{align}

 With reference to the above system: after providing a discussion of the modeling in Section \ref{cantdeflects} (which includes comparisons against other cantilever models appearing in the applied PDE/control literature), we prescribe a functional setup for the problem in Sections \ref{energiessec} and \ref{spacessec}, with definitions of strong and weak solutions in Section \eqref{sols}.

\section{Cantilever Large Deflections}\label{cantdeflects}

Perhaps the most important distinction for large deflection models of cantilevered structures versus fully clamped or hinged structures (completely restricted along the boundary) is that of extensibility. In the case of {\em extensible} beams, transverse deflection necessarily leads to local stretching, which is a principal contributor to the nonlinear elastic restoring force; in the case of clamped-free conditions, for instance, the engineering literature indicates that the beam should be taken to be {\em inextensible} \cite{paidoussis, inext1,semler2}. This includes recent aerodynamic experiments  \cite{dowell4,inext2} suggesting that extensibility (stretching on bending) is not the dominant nonlinear effect in cantilevers.

The property of inextensibility is best characterized as local arc length preservation throughout deflection. Letting $u(x,t)$ and $w(x,t)$ correspond, as before, to the  longitudinal and transverse Lagrangian deflections, the condition manifests itself in the requirement:
$$w_x^2+(1+u_x)^2=1.$$
We note that, if both $w(x=0)$ and $w(x=L)$ are zero, then the beam {\em must be} extensible in order to deflect. 
In the diagrams below, we consider a cantilever in an axial flow, with given unperturbed flow field $\mathbf U = \langle U , 0 \rangle$.

  \begin{center}
\scalebox{1}{
\begin{tikzpicture}[scale=1.2,>=latex']
\draw[->,thick](0,0)--(6,0) node[below]{\scriptsize $x$};
\draw[->,thick](0,0)--(0,2.25) node[left]{\scriptsize $z$};
\draw[gray,line width=3pt] (0,0)--(5,0);
\draw (5,0.25)--(5,-0.25)node[below]{\scriptsize $L$};
\draw (0,0.25)--(0,-0.25)node[below]{\scriptsize $0$};
\draw[line width=3pt,smooth] plot[domain=0:4.5]({\x},{1.25-1.25*cos(0.125*pi*\x r)});
\draw[<->] (4,0)--node[pos=0.5,left]{\scriptsize $w$}(4,1.25);  
\draw[thick,dashed](4.5,0)--(4.5,2.2);
\draw[thick,dashed](5,0)--(5,2.2);
\draw[<->] (4.5,1.5)--node[pos=1,right]{\scriptsize $u(L)$}(5,1.5);
\draw[dashed,->] (0.25,2.1)--(2.5,2.1);  
\draw[dashed,->] (0.25,1.4)--(2.5,1.4);  
\draw (1.5,1.75) node{\scriptsize $\mathbf{U} \equiv \langle U, 0\rangle$};
\end{tikzpicture}

\begin{tikzpicture}[scale=1.2,>=latex']
\draw[->,thick](0,0)--(6,0) node[below]{\scriptsize $x$};
\draw[->,thick](0,0)--(0,2.25) node[left]{\scriptsize $z$};
\draw[gray,line width=3pt] (0,0)--(5,0);
\draw (5,0.25)--(5,-0.25)node[below]{\scriptsize $L$};
\draw (0,0.25)--(0,-0.25)node[below]{\scriptsize $0$};
\draw[line width=3pt,smooth] plot[domain=0:5]({\x},{1.25-1.25*cos(0.125*pi*\x r)});
\draw[<->] (4,0)--node[pos=0.5,left]{\scriptsize $w$}(4,1.25);  

\draw[thick,dashed](5,0)--(5,2.2);

\draw[dashed,->] (0.25,2.1)--(2.5,2.1);  
\draw[dashed,->] (0.25,1.4)--(2.5,1.4);  
\draw (1.5,1.75) node{\scriptsize $\mathbf{U} \equiv \langle U, 0\rangle$};
\end{tikzpicture}}
\end{center}

\subsection{Extensible Cantilevers}

The first model we describe, for context, is a baseline linear model in {\small $w$}---the cantilevered Rayleigh beam \cite{beams}: 
\begin{equation}\label{linearplate}
\begin{cases} (1-\alpha\partial_x^2)w_{tt}+D\partial_x^4w +[k_0-k_1\partial_x^2+k_2\partial_x^4]w_t = p(x,t) \\
w(t=0)=w_0;~~w_t(t=0)=w_1 \\
w(0)=w_x(0)=0;~~w_{xx}(L)=0,\quad\partial_x\big[\alpha w_{tt}-Dw_{xx}+k_1w_t-k_2\partial_x^2w_t\big]\Big|_{x=L}=0. 
\end{cases}
\end{equation}
Above, we have (after a traditional non-dimensionalization \cite{lagleug,beams}) the additional physical quantities:
\begin{itemize}  \setlength\itemsep{.01cm}
\item $\alpha \ge 0$: rotary inertia coefficients in beam filaments; $\alpha=0$ gives the traditional Euler-Bernoulli beam;
\item $k_i \ge 0$ damping coefficients; $k_0$ represents {\em weak} damping, $k_1$ represents {\em square root-like} damping, and $k_2$ (as before) is {\em Kelvin-Voigt} type damping.
\end{itemize}
\begin{remark}
Note that when $\alpha=k_1=k_2=0$, the traditional Euler-Bernoulli cantilever is recovered. If $\alpha=k_1=0$ but $k_2>0$, the third order boundary condition 
$$Dw_{xxx}(L)+k_2w_{xxxt}(L)=0~~\implies~~w_{xxx}(L)=0~~\forall~~t >0,$$ hence Kelvin-Voigt damping can be considered with traditional free end boundary conditions of
$$w_{xx}(L)=w_{xxx}(L)=0.$$ \end{remark}
In the Rayleigh ($\alpha>0$) or Euler-Bernoulli ($\alpha=0$) beams above, there is no evolution for the in-plane displacement $u(x,t)$. 
 The extensible  cantilever system found in \cite{lagleug}, in addition to standard elasticity assumptions, invokes a quadratic strain-displacement law. As a system,  it is nonlinearly coupled in $u$ and $w$. More specifically, the evolutions in $w$ and $u$ employ nonlinear restoring forces resulting from the beam's {\em extension}:
 \begin{equation} \label{LLsystem}
\begin{cases} u_{tt} -D_1\partial_x\left[u_x+\frac{1}{2}(w_x)^2\right]=0 \\ \ds
(1-\alpha\partial_x^2)w_{tt}+D_2\partial_x^4 w -D_1\partial_x\left[ w_x(u_x+\frac{1}{2}w_x^2)\right] = p(x,t) \\
u(0)=0;~~\left[u_x(L)+\frac{1}{2}w_x^2(L)\right]=0 \\
w(0)=w_x(0)=0;\\
w_{xx}(L)=0, ~~-\alpha\partial_xw_t+D_2w_{xxx}(L)=0\\
u(t=0)=u_0;~~u_t(t=0)=u_1;~~w(t=0)=w_0;~~w_t(t=0)=w_1. \end{cases}
\end{equation}
This {\em Lagnese-Leugering} system is the beam analog of the so called full von Karman plate equations \cite{koch}.
Above, $D_1,D_2 >0$ are two independent stiffness parameters, and we have taken the beam without damping effects.

\begin{remark}
 In the unscaled version of the equations $\ds D_1=\dfrac{E}{\rho}$, $\ds \alpha=\dfrac{I}{A}$, and $\ds D_2=\dfrac{EI}{\rho A}$, where $\rho$ is the mass density (per unit volume) of the beam, $I$ is the beam's moment of inertia w.r.t. the $y$-axis, $E$ is the Young's modulus, and $A$ is the cross-sectional area of the beam at rest.  
\end{remark}

One further consideration can be made as a simplification of the above system when we take in-plane accelerations to be negligible, $u_{tt} \approx 0$. Then  we have
$$u(L)-u(0)=c(t)L-\frac{1}{2}\int_0^L w^2_x(\xi)d\xi.$$
We can impose the assumption that the in-plane displacements at the free end of the beam must remain fixed: 
$u(0,t)=0,$ and $u(L,t)=C$, where $C>0$ represents initial in-plane stretching, and $C<0$ compression. As a result, we see that $\ds c =\frac{C}{L}+ \frac{1}{2L}\int_0^L w_x^2(\xi)d\xi.$ Plugging this back into \eqref{LLsystem}, we obtain a scalar extensible cantilever, as studied in \cite{HTW}:
 \begin{equation}\label{Bergerplate}
\begin{cases} (1-\alpha\partial_x^2)w_{tt}+D_2\partial_x^4 w +(k_0-k_1\partial_x^2) w_t-\left[\dfrac{D_1C}{L}+\dfrac{D_1}{2L}\|w_x\|^2\right]w_{xx} = p(x,t)\\
w(t=0)=w_0;~~w_t(t=0)=w_1 \\
w(0)=w_x(0)=0;~~w_{xx}=0; \\
-\alpha\partial_x[w_{tt}+k_1w_t]+D\partial_x^3 w+(b_1-b_2\|w_x\|^2)w_x=0~\text{ at } x=L. \end{cases}
\end{equation}
In the reference \cite{HTW}, a principal component of the analysis is whether $\alpha>0$ or $\alpha=0$. In the case where $\alpha>0$, the results are strong, but stabilization estimates require the damping strength to be tailored to the inertia, i.e., when $\alpha>0$, then $k_1>0$.

\subsection{The Inextensible Cantilever}\label{inextcant}

In each of the models in the previous section, the beam is permitted to be extensible, with nonlinear effects measured in terms of the beam's (local) extension. In contrast, let us now consider {\em inextensibility}. Let $\varepsilon(x)$ denote the {\em axial strain} along the beam centerline. Then, classically, we have the relation \cite{semler2,inext2}:
\begin{equation}
\label{straincons}
[1+\varepsilon(x)]^2=(1+u_x)^2+w_x^2.
\end{equation}

Since we are considering an {\em inextensible} beam, the arc length is preserved throughout deflection and thus we should consider $\varepsilon(x)=0$. Hence the full inextensibility condition can be written as:
\begin{equation}\label{fullinext}
1=(1+u_x)^2+w_x^2.
\end{equation}

Now, let us define the potential energy ($E_P$) with $\varepsilon(x)=0$ via beam curvature $\kappa$ and stiffness $D$ (flexural rigidity) \cite{semler2} in the standard way:
 $$E_P \equiv \frac{D}{2}\int_0^L\kappa^2dx.$$  
 (For a more rigorous derivation of the potential energy, see \cite{newkevin1,follower,langre2}.) The familiar expression for curvature gives, in this instance:
 $$\kappa = \dfrac{(1+u_x)w_{xx}-w_xu_{xx}}{\left((1+u_x)^2+w_x^2\right)^{3/2}}.$$ Invoking the inextensibility constraint\eqref{fullinext}, we obtain:
 $$\kappa=(1+u_x)w_{xx}-w_xu_{xx}.$$
 Via \eqref{fullinext}, we can simplify $\kappa$ to an expression only in $w$: 
 \begin{equation}\label{kappa}
 \kappa=w_{xx}[1-w^2_x]^{-1/2}.
 \end{equation}
 It is at this point we invoke simplifications in both the inextensibility constraint \eqref{fullinext} and the curvature $\kappa$.
 
First, in the inextensibility constraint, we retain the term $w_x^2$, but drop the term $u_x^2$ in \eqref{fullinext}, resulting in the {\em effective inextensibilty constraint}:
\begin{equation}\label{inext}
u_x=-\frac{1}{2}w_x^2.
\end{equation}

 Now, as is standard in elasticity theory, we approximate the curvature $\kappa$ via a Taylor expansion. In line with the above order considerations (for consistency with \eqref{inext} \cite{inext1,semler2}), we retain terms up to order $w_x^2$, yielding: 
 $$ \kappa^2=w^2_{xx}[1-w^2_x]^{-1} \approx w^2_{xx}(1+w_x^2).$$
 \begin{remark} Note that this is the key point which distinguishes various theories of nonlinear elasticity; in linear elasticity $\kappa \approx w_{xx}$.\end{remark}

 With this analysis, the the potential energy becomes:
 $$E_P = \frac{D}{2}\int_0^L w_{xx}^2 \left( 1+ w_{x}^2 \right)dx.$$ 
 The kinetic energy ($E_K$) is defined in the standard way as: 
 $$E_K = \frac{1}{2}\int_0^L\left( u^2_{t} + w^2_{t} \right) dx.$$ 
	
To derive the equations of motion and the associated boundary conditions, Hamilton's Principle is utilized \cite{inext1}. The inextensibility constraint,  $f \equiv u_{x} + (1/2)w^2_{x}=0$, is enforced via a Lagrange multiplier $\lambda$, appended to the system. The Lagrangian is expressed in the usual way:
\begin{equation}
\label{Lagrangian}
\mathcal{L} = E_{K} - E_{P} + \int_0^L \lambda f dx.
\end{equation}
After taking the first variation of $\mathcal L$ and performing the necessary integration by parts with respect to time and space, Hamilton's principle yields the Euler-Lagrange equations of motion and associated boundary conditions \eqref{dowellnon}--\eqref{dowellnon3}. We note that, remarkably, the standard linear clamped-free boundary conditions are obtained. 

\subsection{Discussion of Damping}\label{damping}

 For nonlinear hyperbolic-like problems, the regularity of $w_t$ is at issue; this is especially true for our results here. Additionally, damping provides many useful features, beyond regularization, for non-conservative problems. Discussion of beam damping types goes far back in both the engineering literature \cite{bolotin} and the mathematical literature \cite{russell:93:JMAA,beamdamping}. 
 
Let us refer to \eqref{linearplate}:  weak damping has the form ~{$ k_0 w_t$},  providing no velocity regularization and a damping effect which is uniform in modes \cite{bolotin}. In the elasticity context,  Kelvin-Voigt damping ~{$k_2\partial_x^4 w_t$} is strain-rate type, and mirrors the principal operator; such damping lifts {$w_t \in H^2$}, while transmuting the underlying linear dynamics to be of parabolic type \cite{che-tri:89:PJM}.  Square root-like damping, { $-k_1\partial_x^2w_t$} \cite{fab-han:01:DCDS}, yields {$w_t \in H^1$}, and interpolates between the previous two damping types.
 
 Let us elaborate on square root-like damping: the damping term $\partial_x^2 w_t$ roughly corresponds to half the order of the principal stress operator $\partial_x^4$, if we ignore the  boundary conditions encoded into $\mathcal A$. (We note that for the cantilever configuration, $\mathcal A^{1/2} \neq -\partial_x^2$ \cite{beamdamping}.) The fractional damping concept can be generalized to powers $[\mathcal A]^\theta w_t$ for $\theta\in[0,1]$, which at the two extremes yield the usual weak damping for $\theta=0$ and strong (Kelvin-Voigt/visco-elastic) damping at $\theta=1$.  The square root scenario $\theta=1/2$ for a  system of elastic type was discussed in earlier work \cite{che-tri:89:PJM,beamdamping}, as  this feedback turns out to reproduce energy decay rates empirically observed in elastodynamics. Fractional damping was investigated in \cite{che-tri:89:PJM} for the abstract system 
$$w_{tt}+\mathcal A^{\theta}w_t+\mathcal Aw=0,$$ demonstrating, in particular, that the ensuing evolution semigroup is  analytic if and only if $\theta \geq \frac{1}{2}$.  
 
 We note that square root damping corresponds to modal damping models \cite{dowell}, as one finds frequently in the engineering literature \cite{bolotin,McHughIFASD2019,follower}.  However, the boundary conditions for a given problem affect the physical interpretation of fractional damping for certain values of $\theta$, and in \cite{beamdamping} it is noted that $\theta=1/2$ has a questionable physical interpretation for a cantilevered configuration. The square root-like damping $\partial_x^2 w_t$ also arises naturally in other beam models. Consider, for instance,  the Mead-Markus
model \cite{mea-mar:69:JSV} for a  sandwich beam \cite{fab-han:01:DCDS}.  
 
In what follows below, the presence of damping---its strength and affect on the regularity of $w_t$---will be critical to the main theoretical results concerning existence and uniqueness for the inextensible dynamics. This is discussed further in the results Section \ref{resultsdiscuss} and Section \ref{optimaldamp}.

\section{PDE Model Studied Here}

With the derivation mentioned above (following \cite{inext1}), we recall the equations of motion, allowing for Kelvin-Voigt damping $k_2 \ge 0$:
\begin{equation}\label{dowellnon*}
\begin{cases} \displaystyle w_{tt}+D\partial_x^4w+ k_2 \partial_x^4 w_t+\mathbf A (w) =p(x,t)& \text{ in } (0,L) \times (0,T)
\\ w(t=0)=w_0(x), w_t(t=0)=w_1(x) \\ w(x=0)=w_x(x=0)=0;~w_{xx}(x=L)=w_{xxx}(x=L)=0. 
\end{cases}
\end{equation}
with the nonlinear, nonlocal $\mathbf A$ given through
\begin{align}\label{dowellnon2*}
\mathbf A (w) =&- \sigma D\partial_x\big[(w_{xx})^2w_x \big]+ \sigma D\partial_{xx}\big[w_{xx}\big(w_x^2\big)\big] + \iota \partial_x\left[w_x\int_x^L u_{tt}(\xi)d\xi \right] \\
u(x)=&~-\frac{1}{2}\int_0^x \left [w_x(\xi) \right ]^2d\xi.\label{dowellnon3*}
\end{align}
To simplify terminology, we use the following language  from here on:
\begin{align*}\text{\bf  \small [NL Stiffness]}=&~ -D\partial_x\big[(w_{xx})^2w_x \big]+D\partial_{xx}\big[(w_x)^2w_{xx}\big] \end{align*} as the {\em nonlinear stiffness terms}, while we refer to 
\begin{align*}\text{\bf  \small [NL Inertia]}=&~ \partial_x\left[w_x\int_x^L u_{tt}(\xi)d\xi \right ]\end{align*} as the {\em nonlinear inertial term} (which is nonlocal, when written in $w$).
We have introduced flags, $\iota,\sigma=0 \text{ or } 1$, in \eqref{dowellnon2*}, in order to easily turn particular nonlinear effects on or off. This is to say, when $\iota=0$, we say that {\bf \small [NL Inertia]} is turned off.

\begin{remark}
As shown in detail in \cite{inext1}, the {\em standard linear} clamped-free (cantilever) boundary conditions are recovered in the variational procedure. This is somewhat remarkable, noting that we have allowed for broad nonlinear and nonlocal effects. We also point out that when boundary forces are enacted at $x=L$, not only are the free end boundary conditions altered for $w$ and $u$, but the relationship between $w$ and $u$ in \eqref{dowellnon3} is impacted through the Lagrange multiplier $\lambda$ in \eqref{Lagrangian} (itself having boundary conditions). See \cite{follower}.
\end{remark}

\begin{remark}\label{quasi}
For convenience, we note two expansions.
First:
$$\text{\bf \small [NL ~Stiffness]}=-D\partial_x\big[(w_{xx})^2w_x \big]+D\partial_{xx}\big[(w_x)^2w_{xx}\big]=D[w_{xxx}^3+4w_xw_{xx}w_{xxx}+w_x^2w_{xxxx}],$$
which highlights the quasilinear nature of the PDE (with high order semilinearity).

Second:
$$\text{\bf \small [NL ~Inertia]}= \partial_x\Big[w_x\int_x^Lu_{tt}d\xi\Big] =w_{xx}\int_x^L u_{tt}d\xi-w_xu_{tt},~~\text{with}~~ u_{tt}=-\int_0^x[w_{xt}^2+w_xw_{xtt}]d\xi.$$
This expansion highlights the fact that the system can be closed in $w$, as well as the high velocity regularity required to interpret the strong form of the PDE. 
\end{remark}

\section{Theory of Solutions}\label{theory}
\subsection{Notation}

For a given domain $D$,
its associated $L^{2}(D)$ will be denoted as $||\cdot ||_D$ (or simply $||\cdot||$ when the context is clear). Inner products in a Hilbert space  are written $(\cdot ,\cdot)_{H}$ (or simply $(\cdot ,\cdot)$ when $H=L^2(D)$ and the context is clear). We will also denote pertinent duality pairings as $\left\langle \cdot ,\cdot \right\rangle _{X\times X^{\prime }}$, for a given
Hilbert space $X$. The space $H^{s}(D)$ will denote the standard Sobolev space of
order $s$, defined on a domain $D$, and $H_{0}^{s}(D)$ denotes the closure
of $C_{0}^{\infty }(D)$ in the $H^{s}(D)$-norm $\Vert \cdot \Vert
_{H^{s}(D)}$, also written as $\Vert \cdot \Vert _{s}$. For $\Gamma \subset \partial D$, boundary restrictions $u\big|_{\Gamma}$ are taken in the sense of the trace theorem for $u \in H^{1/2+}(D).$
\subsection{Energies}\label{energiessec}
With reference to Section \ref{inextcant}, we have the following energies:
\begin{equation}
E(t) \equiv E_K(t)+E_P(t) \equiv \frac{1}{2}\left[||w_t||^2+\iota ||u_t||^2\right]+\dfrac{D}{2}\left[||w_{xx}||^2+\sigma ||w_xw_{xx}||^2\right].
\end{equation}
Note that the energies include the flags.

As we suppress the $u$ variable here
via the effective inextensibility constraint \eqref{inext},
we can write the energetic terms in $w$ along, with separated by linear and nonlinear designations:
\begin{equation}
E(t) = E_L(t)+E_N(t) \equiv \dfrac{1}{2}||w_t||^2+\dfrac{D}{2}||w_{xx}||^2+\dfrac{\sigma D}{2}||w_xw_{xx}||^2+\dfrac{\iota}{2}\left|\left|\int_0^xw_{x}(\xi)w_{xt}(\xi)d\xi\right|\right|^2.
\end{equation}
In the unforced situation with $p(x,t) \equiv 0$, we note that the formal energy identity is obtained by the velocity multiplier $w_t$ on \eqref{dowellnon} and 
$$E(t) +k_2\int_s^t||w_{xxt}||_{L^2(0,L)}^2d\tau= E(s).$$

In what follows we will define higher order energies corresponding to the topology of smooth solutions. This will be done in the corresponding sections.

\subsection{Spaces}\label{spacessec}
The principal displacement state space for cantilevered beam dynamics takes into account the clamped conditions:
\[
  H^2_* = \{ v \in H^2(0,L) :  v(0) =0,\quad v_x(0) = 0 \}.
\]
This space is equipped with an $H^2$ equivalent inner product
\begin{equation}\label{H2*prod}
   (v,w)_{H^2_*} = D (v_{xx}, w_{xx}).
 \end{equation}
Denoting $R$ as the Riesz isomorphism $H^2_*\to [H^2_*]'$, it is given by:
 \begin{equation}\label{def:R}
 R(v)(w) \equiv  (v,w)_{H^2_*}\,.
 \end{equation}
This framework is conveniently induced by the generator of the linear cantilever dynamics:
\[
  \cA : \mathcal D(\cA)\subset L^2(0,L) \to L^2(0,L),~~ \cA f \equiv D\partial_x^4 f,  
\]
\begin{equation}\label{def:cA0}
\mathcal D(\cA) = \{ w \in H^4(0,L) : w(0) =w_x(0)=0, \; w_{xx}(L)=w_{xxx}(L) = 0\}.
\end{equation}
From this we have in a standard fashion: 
\[
  \mathcal D(\cA^{1/2}) = H^2_*, \quad \mathcal D(\cA^{-1/2}) = [H^2_*]' ~\text{ and }~ \cA^{1/2} = R \quad \text{(as in \eqref{def:R})}.
\]
Then $(u,\cdot)_{H^2_*}$ is the extension of $(\cA u, \cdot)$ from $\mathcal D(\cA)$ to $H^2_*$ which gives \eqref{H2*prod}. 

\begin{remark} We note that, despite the topological equivalence of $\mathcal D(\mathcal A^{1/2})$ and $H^2_*$, it is not the case that $\mathcal A^{1/2}$ can be identified with $-\partial_x^2$ on $H^2_*$ \cite{beamdamping}. This relates to a deep discussion connected to {\em square root-like damping}, as described below in Section \ref{damping}. \end{remark}

Using the above spaces we can define the appropriate state space(s) for our dynamics. The finite energy space will be denoted as: $$\mathscr H \equiv H^2_*\times L^2(0,L),$$ with the inner product: $y=(y_1,y_2),~ \tilde y=(\tilde y_1,\tilde y_2) \in \mathscr H$
\begin{equation}\label{Hal-prod}
  (y,\tilde y)_{\mathscr H} = (y_1,\tilde y_1)_{H^2_*} + (y_2,\tilde y_2)_{L^2}.
\end{equation}
We note that the norm in $\mathscr H$ topologically corresponds to the energy functional $$E_L(t) = \dfrac{D}{2}||w_{xx}||^2+\dfrac{1}{2}||w_t||^2.$$

In our discussions, we will also require a stronger state space (corresponding to {\em strong} solutions):
\begin{equation}\label{strongspace}
\mathscr H^s \equiv \begin{cases} \mathcal D (\mathcal A) \times \mathcal D(\mathcal A^{1/2}), & \iota=k_2=0, \\[.2cm] 
\mathcal D (\mathcal A) \times \mathcal D(\mathcal A), & \iota=1, k_2>0.
\end{cases}
\end{equation}
The norm in $\mathcal H^s$ is taken (equivalent to the natural operator-induced norm) to be:
$$||y||_{\mathcal H^s}^2 =
\begin{cases}  ||\partial_x^4 y_1||_{L^2}^2+||\partial_{xx}y_2||_{L^2}^2, & \iota=k_2=0, \\[.2cm] 
 ||\partial_x^4 y_1||_{L^2}^2+||\partial_{x}^4y_2||_{L^2}^2,& \iota=1, k_2>0.
\end{cases}
$$

\subsection{Definition of Solutions}\label{sols}
We provide the natural setting for the weak formulation of the problem; this will yield the appropriate starting point for our numerical (modal) methods, as well as provide the appropriate abstract setting for analysis of the equations of motion. Ultimately, we will construct weak solutions that possess additional regularity; these will turn out to be strong solutions.

Then, the associated weak form of \eqref{dowellnon} has the form:
\begin{align} \label{weakform}
\dfrac{d}{dt}\Big[(w_t,\phi)&+ \iota \Big(\int_0^xw_xw_{xt}d\xi,\int_0^xw_x\phi_xd\xi \Big)\Big]- \iota \Big(\int_0^xw_xw_{xt}d\xi,\int_0^xw_{xt}\phi_xd\xi\Big) \\
&+k_2(w_{xxt},\phi_{xx})+D(w_{xx},\phi_{xx})+ \sigma D(w_{xx}w_x,w_{xx}\phi_x)+ \sigma D\big(w_{xx}w_x,\phi_{xx}w_x\big) = (p,\phi),\nonumber
\end{align}
for $\phi \in H^2*$, and where the $d/dt$ above is interpreted in the sense of $\mathscr D'(0,T)$. 
When $\sigma>0$, the {\bf \small [NL Stiffness]} is in force; similarly, when $\iota>0$, {\bf \small [NL Inertia]} is in force.
When $k_2>0$, Kelvin-Voigt damping is imposed.

We now give precise definitions of solutions making reference to the weak form \eqref{weakform} above:

\begin{definition}\label{sol1}
We say a {\em weak} solution to \eqref{dowellnon}, with $k_2 = \iota = 0$ and $ \sigma=1$ is a function $w$, with
\begin{equation*}
w \in L^2\left(0,T;H_*^2(0,L)\right);~w_t \in L^2\left(0,T;L^2(0,L)\right);~w_{tt} \in L^2\left(0,T;[H_*^2(0,L)]'\right)
\end{equation*}
that satisfies \eqref{weakform}.

Moreover, for any $\chi \in H^2_*$, $\psi\in L^2(0,L)$, we require
  \begin{equation}\label{weak-ic}
    (w,\chi)_{H^2_*}\big|_{t\to 0^+} = (w_0, \chi)_{H^2_*},\quad (w_t,\psi)\big|_{t\to 0^+} = (w_1,\psi).
  \end{equation}

\end{definition}

\begin{definition}\label{sol2}
A {\em weak} solution to \eqref{dowellnon} with $k_2>0$ and $\iota=\sigma=1$ is a function $w$, with
\begin{equation*}
w \in L^2\left(0,T;H_*^2(0,L)\right); ~ w_t \in L^2\left(0,T;H^2_*(0,L)\right); ~ w_{tt} \in L^2\left(0,T;[H_*^2(0,L)]'\right),
\end{equation*}
such that \eqref{weakform} holds.

Moreover, for any $\chi \in H^2_*$, $\psi\in L^2(0,L)$, we require
  \begin{equation}\label{weak-ic*}
    (w,\chi)_{H^2_*}\big|_{t\to 0^+} = (w_0, \chi)_{H^2_*},\quad (w_t,\psi)\big|_{t\to 0^+} = (w_1,\psi).
  \end{equation}
\end{definition}

\begin{remark}
For $k_2>0$ and $\iota>0$, the definition of weak solution is self-consistent; this is to say, for such a function $w$, all terms in \eqref{weakform} are well-defined. 
We note that for $k_2=0$, there are issues with the a priori regularity of $w_t \in L^2(0,T;L^2(0,L))$ and the interpretation of the {\bf \small [NL Inertia]} terms. \end{remark}

Now, for strong solutions:
\begin{definition}
A {\em strong} solution to \eqref{dowellnon} with $k_2 = \iota = 0$ and $ \sigma=1$ is a weak solution (as in Definition \ref{sol1}) with the additional regularity
$$w \in L^2 \left(0,T; \mathcal D(\mathcal A)\right);~w_t \in L^2(0,T;H^2_*(0,L));~w_{tt} \in L^2 \left( 0,T; L^2(0,L) \right).$$ 
\end{definition}

\begin{definition}
A {\em strong} solution to \eqref{dowellnon} with $k_2>0$, $\iota=\sigma=1$ is a weak solution (as in Definition \ref{sol2}) with the additional regularity
$$w \in L^2 \left(0,T;\mathcal D(\mathcal A) \right);~w_{t} \in L^2 \left( 0,T; \mathcal D(\mathcal A))\right);~w_{tt} \in L^2 \left( 0,T; H^2_{*}  \right).$$
\end{definition}

\subsection{Well-posedness Results}\label{wellpresults}

In this section we state recent theoretical results about {\em strong} solutions to \eqref{dowellnon}. The proofs of these theorems appear in \cite{maria}, with an effort to have a streamlined presentation of the underlying modeling and theory supporting the numerical simulations below.

We begin with a simple well-posedness result for the nonlinear stiffness portion of the model.
\begin{theorem}\label{withoutiota}
Take $\sigma =1$ with $\iota=k=0$, and consider $p,~p_t,~p_{xx} \in L^2(0,T;L^2(0,L))$. For smooth data $w_0 \in \mathcal D(\mathcal A)$, $w_1 \in H_*^2$, strong solutions exist up to some time $T^*(w_0,w_1)$. For all $t \in [0,T^*)$, the solution $w$ is unique and obeys the energy identity:
$$E(t) = E(0)+\int_0^t (p,w_t)_{L^2(0,L)}d\tau.$$ 
Solutions depend continuously on the data in the sense of $C([0,T^*);\mathscr H)$ with an estimate on the difference of two trajectories, $z=w^1-w^2$:
$$\sup_{t \in [0,T]}\big|\big|(z(t),z_t(t))\big|\big|_{\mathscr H} \le C\Big(||(w^i(0),w^i_t(0))||_{\mathscr H^s},T\Big)\big|\big|\big(z(0),z_t(0)\big)\big|\big|_{\mathscr H},~~\forall~T \in [0,T^*).$$
\end{theorem}

Now, we consider the entire model---$\iota=1$ and $k_2>0$. 
\begin{theorem}\label{withiota}
Take $\sigma=\iota=1$ and $k_2>0$, and consider $p,~p_t,~p_{xx} \in L^2(0,T;L^2(0,L))$. For initial data $w_0, w_1  \in \mathcal D(\mathcal A)$, strong solutions exist up to some time $T^*(w_0,w_1)$. For all $t \in [0,T^*)$, the solution obeys the energy identity:
$$E(t)+k_2\int_0^t||w_{txx}||^2_{L^2(0,L)} = E(0)+\int_0^t (p,w_t)_{L^2(0,L)} d\tau.$$
\end{theorem}
Note that we make no claims of uniqueness above when $\iota=1$. 
\begin{remark}
The dependence $t^*(w_0,w_1)=t^*\big(||(w_0,w_1)||_{\mathcal H^s}\big)$, indicates that the time of existence depends on the size of the initial data in the strong norm $\mathcal H^s$. The definition of this space depends on the values of $\iota, k_2$; see \eqref{strongspace}. 
\end{remark}

\subsection{Previous Results and Discussion}\label{resultsdiscuss}

In Section \ref{cantdeflects} above, we provided a discussion of the modeling and previous mathematical analyses of cantilevers and cantilever large deflections. Section \ref{inextcant} provides a discussion of inextensibility, along with early references towards its modeling and recent engineering-oriented numerical references. Finally, Section \ref{damping} provides a discussion of damping mechanisms in beams, with a focus on cantilevers. Now, in this section, we provide a brief discussion of the remaining relevant literature from the point of view of the results presented above.

      The earliest modeling and computational work concerning inextensibility seems to come from Pa\"{i}doussis et al. \cite{paidoussis,semler2} in the context of pipes conveying fluid.  Regularizing higher order  Kelvin-Voigt $(k_2>0)$ damping was included in these structural models when producing numerical results.  The recent paper \cite{inext1} provides the Lagrange multiplier derivation of the inextensible beam \eqref{dowellnon*}--\eqref{dowellnon3*} (which we follow in our modeling discussions here); \cite{inext1} also produces the so called Rayleigh-Ritz modal equations of motion, which are studied as a nonlinear ODE system. These approaches are considered and modified in the presence of non-conservative forces in latter papers, such as the nonlinear piston theory \cite{McHughIFASD2019} or non-conservative follower forces \cite{follower,langre2}.
      Apart from their modeling aspects, these papers are primarily numerical in nature, focusing on the onset and qualitative properties of dynamic instability. 
      The very recent \cite{newkevin1} provides a thorough review of modern nonlinear beam theories, taken from the engineering point of view (of a similar ilk as the earlier \cite{beams}).
      The engineering literature shows that the inextensible model described in this paper  performs well \cite{dowell4,inext2,paidoussis3,experimental} when compared to experiments.

Moving on to the mathematical literature, we assert that---to the knowledge of the authors---{\em no PDE or control-theoretical work has treated the inextensible beam}. This is to say, there seems to be no available existence and uniqueness theory for \eqref{dowellnon*}--\eqref{dowellnon3*}. We attempt to remedy this with our results above.  Our results in Section \ref{wellpresults} provide local well-posedness for the case of nonlinear stiffness only ($\sigma=1,~\iota=0$), without the need for damping. On the other hand, to obtain a local well-posedness result for the case where nonlinear inertia is present ($\sigma=1,~\iota=1$), we must add strong (Kelvin-Voigt type) damping, $k_2>0$, and adjust the state-space accordingly \cite{redbook}. 

First, we note that our well-posedness results are consistent with what is to be expected for quasilinear beams and plates (see, e.g., \cite{ong,xiang}), owing to the fact that the {\bf \small [NL Stiffness]} is quaslinear in nature. Moreover, the uniqueness in Theorem \ref{withoutiota} is nontrivial, based on exploiting the particular polynomial structure of the {\bf \small [NL Stiffness]} term. 
Although these results may seem weaker than one might expect, we note that these represent the first such existence and uniqueness type results for the inextensible cantilever. Moreover, it is clear (and explained below) that without damping of the form $\mathcal A^{\theta}w_t$ (for $\theta>0$ sufficiently large), there is no hope of closing energy estimates when $\iota=1$. So although the need for damping $k_2>0$ in our second existence result seems odd, it is precisely because of the nonlocal quasilinear nature of {\bf \small [NL Inertia]} terms. Indeed the effective inextensibility constraint provides a relation between $u$ and $w$, and when this is expanded---see Remark \ref{quasi}---it is clear that higher regularity of $w_t$ is necessary to interpret the solution. This can also be seen from the weak form \eqref{weakform}. The {\bf \small [NL Inertia]} term prevents the equations of motion in \eqref{dowellnon*}--\eqref{dowellnon3*} from being written as a traditional second order evolution in time (for $w$), and truly distinguish this model from previous beam models.

      It is worth emphasizing that, although there is extensive literature for nonlinear beams (and plates), to the best of the authors' knowledge, there is very little mathematical discussion of {\em nonlinear} cantilevered beams at all. One can consult the aforementioned paper \cite{lagleug} on semigroup well-posedness and stabilizability of extensible beam systems, along with the related (simplified) models in \cite{HTW} and \cite{beam4}.
  
 Numerically, we utilize a dynamic modal approach, akin to what is standard in the aeroelasticity literature (for instance, \cite{vedeneev}). Via this approach, a fully nonlinear (implicit) system of ODEs is obtained by expanding the solution in in-vacuo mode shapes and implementing a Galerkin procedure to determine time-dependent Fourier coefficients. Owing to the complex nature of the nonlinearities for the inextensible beam model, finite difference methods are not developed here, as are used, for instance, in the beam study \cite{HHWW}, which compares modal methods and spatially discretized methods.

\subsection{A Priori Estimates and Comments on Well-Posedness Proofs}
In this section we remark briefly on the a priori estimates associated with \eqref{dowellnon*}--\eqref{dowellnon3*}, and the corresponding construction of solutions. Full details appear in the mathematically-oriented paper \cite{maria}; we suffice here to provide an overview of the well-posedness strategy and accompanying scheme for construction of solutions and energy estimates.

\subsubsection{Stiffness Only}
The strategy we follow for obtaining well-posedness, is to firstly establish existence for the quasilinear {\bf \small [NL Stiffness]} component ($\sigma=1$, $\iota=k_2=0$). Following a standard tack, we utilize a Galerkin procedure, taking the standard spatial Fourier basis for the linear, in vacuo beam dynamics on $H^2_*$. Upon implementing the Galerkin procedure, we obtain approximate solutions that satisfy the finite dimensional analog of \eqref{weakform}. The baseline energy identity at the finite energy level $H^2_* \times L^2(0,L)$ yields associated weak limit points. The 1-D Sobolev embedding for {$H^2_*$} provides {$w_x \in L^{\infty}$}, which is adequate to identify weak limits for the term {$\Big([w_x^n]^2w^n_{xx},\phi_{xx}\Big)$} in \eqref{weakform}, with $\phi \in H^2_*$; however, the term {$\Big([w_{xx}^n]^2w^n_{x},\phi_{x}\Big)$} is more delicate. One obtains a limiting measure as the $*$-weak limit point via the Alaoglu Theorem, but additional compactness is needed to identify it in {$H_*^2$} and associate it with a weak solution. (In this setting, direct use of the Dunford-Pettis criterion is not amenable.) Hence, to obtain the needed compactness, we work  with {\em smooth solutions}, in line with standard quaslinear theory. 

Specifically, we employ energy methods for higher order (differentiated) equations and exploit the polynomial structure of the nonlinear terms associated with {\bf \small [NL Stiffness]}. Careful use of a sequence of multipliers, along with delicate estimation of nonlinear terms via interpolation and Sobolev theorems, yields the necessary energy estimates which we now describe.

Let us define $$E_0(t) = \frac{1}{2} ||w_{t}||^2 + \frac{D}{2}||w_{xx}||^2 + \frac{D}{2}||w_xw_{xx}||^2,$$ corresponding to the identity obtained by the velocity multiplier $w_t$ in the equations \eqref{dowellnon*}--\eqref{dowellnon3*} with $\iota=k_2=0$. We have the  corresponding estimate immediately:
\begin{equation}
\label{firstlevel}
E_0(t) = E_0(0) + \int_0^t  \left(p , w_{t} \right) d\tau ~~\text{for all}~~ t>0.
\end{equation}
\begin{remark}
Note that the above estimate is the same one presented in \ref{energiessec}, omitting the nonlinear inertial part. \end{remark}

Now, letting $$E_{1}(t) = ||w_{tt}||^2 +  ||w_{xxt}||^2 + ||w_{xx}w_{xt}||^2 +   ||w_{xxt}w_{x}||^2$$ be the energy corresponding to the {\em time} differentiated version of the stiffness-only equation $(\iota=0)$, and taking the $w_{tt}$ multiplier, we obtain:
\begin{equation}
\label{timediff}
E_1(t) \leq f_1\Big(p_t,E_0(0), E_1(0)\Big) + f_2\Big(p,E_0(0)\Big)t + c \int_0^t E^2_1(\tau)d\tau,
\end{equation}
where $c>0$ and the $f_i$ are smooth, real-valued functions of their arguments. By dependence on $p$ we mean dependence on the norm
$\displaystyle ||p||_{L^2(0,t; L^2(0,L))}$ (mutatis mutandis for derivatives of $p$, such as $p_t,p_{xx}$).

Using a standard {\em nonlinear} version of Gr\"{o}nwall's lemma \cite{gronwall} we obtain a local-in-time estimate:
\begin{equation}
\label{NonlinGronw}
E_1(t) \leq  \frac{ f_1  + f_2t }{1- c \left[ f_1 t + f_2 t^2 \right]}~~0 \leq t \leq T^* ~~\text{where}~~ T^* = \sup_{t} \left \{ c \left[ f_1 t + f_2 t^2 \right] <1 \right \}. 
\end{equation}

From \eqref{NonlinGronw} we can deduce that the Galerkin approximations satisfy a local-in-time bound, providing boundedness in the associated norms of $E_0$ and $E_1$ for some finite time depending on the initial data\footnote{Or, conversely, given any time $T$, there is a ball of initial data sufficiently small in the sense of $E_i(0)$ for which solutions exist up to $T$.}. Unlike standard semilinear theory, we cannot obtain the needed regularity on $\partial_x^4w$ through the equation with the additional regularity of $w_{tt} \in L^{\infty}(0,T^*;L^2(0,L))$. To obtain the final a priori bound for the {\bf \small [NL Stiffness]}, we define $$V(t) =  ||\partial_x^4w||^2 + ||w_{x}\partial_x^4||^2 + ||w_{xx}w_{xxx}||^2,$$ corresponding to the conserved quantity associated to {\em two}  {\em space} differentiations, taken with the $w_{xx}$ multiplier. This yields the inequality:
\begin{equation}
\label{spacediff}
\int_0^t V(s) ds \leq f\big(E_0(0), E_1(0),  E_1(t), p_{xx}\big),
\end{equation}
where $f$ here is increasing in its arguments.

Two space derivatives are utilized, as they constitute a convenient fractional power of $\mathcal A^{1/2}$; we observe that the energy identities associated with one space differentiation result in problematic trace terms that cannot be controlled by the conservative energetic terms. Moreover, \eqref{spacediff} highlights the necessity of firstly having a closed estimate for higher time derivatives of the solution.

The combination of \eqref{firstlevel}, \eqref{NonlinGronw} and \eqref{spacediff} yields the final energy estimate for boundedness of $$||w||_{L^2(0,T^*;\mathcal D(\mathcal A))};~~||w_{tt}||_{L^{\infty}(0,T^*;L^2(0,L))},$$ among others, in terms of initial data ~$E_0(0),~E_1(0), V(0)$.  With additional compactness coming from smooth data, we obtain weak solution with the appropriate limit point identification. Subsequently, with our higher order a priori estimates, we utilize the regularity of the solution to infer that the weak solution is in fact strong. It is an exercise to show that the strong solution satisfies \eqref{dowellnon*}--\eqref{dowellnon3*}  (with $\sigma=1$, $\iota=0$) in a point-wise sense. 

Uniqueness is a nontrivial issue here (of course related to the aforementioned problem of limit point identification). However, we can exploit the polynomial structure of the {\bf \small [NL Stiffness]} terms to obtain a continuous dependence estimate on the initial data, so long as the previous energy estimates hold; from this, uniqueness  follows. To that end, we define $$\mathcal{N}(w)=\partial^2_{x}\left(w_{x}^2 w_{xx} \right) - \partial_x \left(w^2_{xx}w_{x} \right).$$ Let $w^1$ and $w^2$ be two strong solutions of the problem \eqref{dowellnon*} on $t \in [0,T^*)$ with $\sigma=1$, $\iota=k_2=0$ and $z=w^1-w^2$. Then, decomposing the energy multiplier, as applied to the nonlinear difference, we obtain through simple but non-obvious algebraic manipulations:
\begin{align*}
\left( \mathcal{N}(w) - \mathcal{N}(v), z_{t} \right) =&~ \frac{1}{2} \left [ \left(w^2_{xx}, z^2_{x} \right) + \left(w^2_{x}, z^2_{xx} \right) \right ] - \left(w_{xx}w_{xxt}, z^2_{x} \right) -\left(w_{x}w_{xt}, z^2_{xx} \right) \\ &+ \left( v_{x}\left \{ w_{xx}+v_{xx} \right \}, z_{xx}z_{xt} \right) +\left( v_{xx}\left \{ w_{x}+v_{x} \right \}, z_{x}z_{xxt} \right).
\end{align*}
Estimating the above inner-products can be done directly with the help of Cauchy-Schwarz and Young's inequality, invoking the earlier a priori estimates on individual trajectories. This results in a nice energy estimate on $z$ (using the multiplier $z_t$) of the form of the standard (linear-type) Gr\"onwall inequality. Proceeding as is standard, yields continuous dependence in the finite energy topology $\mathscr H$ of the trajectories (and associated uniqueness) for smooth solutions emanating from data in $\mathscr H^s$.

\begin{remark} Note that no damping was needed to obtain uniqueness here. \end{remark}

\subsubsection{With Inertia}\label{w/inertia}
Treating the {\bf \small [NL Inertia]} term perturbatively is challenging, owing to the presence of the term $w_{xtt}$ in the $w$-expanded form of $u_{tt}$---see Remark \ref{quasi}. Hence, we proceed to estimate this term, making use of velocity smoothing associated to the presence of strong (Kelvin-Voigt) type damping with $k_2>0$. We again utilize a Galerkin procedure, though now taking $\sigma=1$, $\iota=1$, and $k_2>0$. 

As one immediately sees from \eqref{weakform}, the weak form peels a time derivative off of all inertial terms; for weak solutions in this situation, we see ``standard" beam requirements for the functions $w_t, ~w_{tt}, ~w_{xx}$ (\cite{HTW} and references therein). As shown in \ref{energiessec},  {$||u_t||_{L^2}$} is formally conserved in the kinetic energy $E_K$, and thus the inextensibility condition provides control of the quantity $$\int_0^xw_xw_{xt}d\xi \in L^{\infty}(0,T;L^2(0,L)).$$  Yet the weak form \eqref{weakform} makes clear that {\em some additional regularity of ~{$w_t$}} is required for appropriately interpreting \eqref{weakform}, namely, so {\small $w_{xt}$}  is well-defined. 

Moreover, the strategy utilized to close estimates for stiffness calls for a time differentiation of the equations, followed by an application of $\mathcal A^{1/2}$ (and associated estimation). For inertial terms, we need to verify that the additional terms are compatible with the aforementioned estimation procedures. Thus, following the prescribed scheme, we cannot ``avoid" differentiating {\bf \small [NL Inertia]} terms in time. This however, showcases the lack of necessary $w_t$ regularity again: 
$$u_{tt}(x)=-\int_0^x \left [w_{xt}^2+w_xw_{xtt} \right ]d\xi,$$
whence we can already see control over the term ~$w_{xtt}$ is necessary even at the undifferentiated equations. To achieve this control, one can attempt to differentiate further in time, but this is ineffective, since every {\em time} differentiation boosts the requisite time regularity of $w_{x}$. Hence, differentiation in time will not provide closed estimates for inertial terms. Repeated spatial differentiation is incompatible with closing the estimates from the earlier {\bf \small [NL Stiffness]}.

Hence, owing to the above discussion, regularity for $w_t$ must be ``borrowed" from some other term in the equation. Some standard regularizations to resolve these sort of issues include the use of linear (Rayleigh-type) rotational inertia ~$w_{tt} \mapsto (1-\alpha\partial_{xx})w_{tt}$ (as in \eqref{linearplate}), which is not helpful for the inextensible model, owing to incompatibility between the {\bf \small [NL Stiffness]} terms and the cantilever boundary conditions.  One might consider utilizing square root-like damping of the form $-k_1\partial_{x}^2w_{t}$ (as in \eqref{linearplate} and \eqref{Bergerplate}), yet for cantilevers, this requires modifying the higher-order boundary conditions; there is some discussion of the physical interpretation of damping mechanisms weak $k_0>0$, square root-like $k_1>0$, and strong $k_2>0$ in \cite{beamdamping}. Here, we proceed with the addition of linear {\em Kelvin--Voigt } damping by taking $k_2>0$ \cite{che-tri:89:PJM,beamdamping}. This is a physically viable form of damping for cantilevers, and it is also used in the engineering literature \cite{semler2}. This choice does not require modification of the higher order beam boundary conditions. 

\begin{remark} It may be the case that a weaker form of $\mathcal A^{\theta}w_t$ damping is sufficient to obtain estimates; we discuss this later in Section \ref{future}. \end{remark}

With the inclusion of strong damping, we may run the procedures: 
$$
\Big\{\partial_t~~/~~\times w_{tt}~~/~~\int_0^L\Big\}~~~~;~~~~
\Big\{\partial_{xx}~~/~~\times w_{xxt}~~/~~\int_0^L\Big\}$$ on \eqref{dowellnon*}
to obtain two additional energy estimates, which close, thanks to $k_2>0$. In this case, the final a priori estimate becomes:
\begin{equation}
\label{finalest}
\mathcal{E}(t) + \mathcal{I}(t) + \mathcal D_0^t[w_t] \leq f_1\big(\mathcal{E}(0),\mathcal{I}(0)\big) +f_2(E_0(0)) \int_0^t \left( \mathcal{E}(s) + \mathcal{I}(s) \right)^2 ds
\end{equation}
for smooth functions $f_i$, increasing in their arguments, with:
$$\mathcal{E}(t) = E_0(t)+E_1(t)+V(t),~~~~~\mathcal{I}(t) = ||u_{t}||^2 + ||u_{tt}||^2 + ||u_{xxt}||^2+\left|\left|\int_0^xw_{xt}^2d\xi\right|\right|^2,$$
and
$$\mathcal D_0^t[w_t]=c(k_2) \int_0^t \left[ ||w_{xxt}||^2 + ||w_{xxtt}||^2 + ||\partial_x^4w_{t}||^2 \right] d\tau.$$
(For clarity, we have above suppressed the dependence on the RHS forcing function $p$ and its derivatives.)
This estimate should be contrasted with the estimates for {\bf [NL Stiffness]} only in \eqref{NonlinGronw} and \eqref{spacediff}, which do not depend on the presence of damping. 
Note that with the addition of damping we obtain a better estimate for $V(t)$; namely, the ``stability-type" multiplier $w_{xx}$ is replaced by the ``energetic" multiplier $w_{xxt}$, which is permitted owing to the presence of damping. Then, similar to \eqref{timediff}, a nonlinear version of Gr\"{o}nwall is utilized as well, providing the final a priori estimates for $\iota=1,~k_2>0$.

Once a priori estimates are in hand, limit passage obtains as before and identification of limits is as in the previous stiffness-only case. We can again use the solution regularity (corresponding to a priori estimates and data requirements in $\mathscr H^s$) to infer that the weak solution is in fact strong, and satisfies the full PDE (with $\sigma=1$ and $\iota=1$) in a point-wise sense. We highlight here that no claims for uniqueness are made for the {\bf \small [NL Inertia]} case $\iota=1$.

\section{Simulation of Inextensible Cantilever Dynamics}
This section is devoted to numerically simulating the inextensible cantilever dynamics. 
In Sections \ref{piston} and \ref{modal} we describe the method and approach to producing dynamic (modal) simulations. Subsequently, in Section \ref{numerics}, we show our numerical results and provide a detailed discussion. Finally, in Section \ref{numconc}, we provide an overview of numerical conclusions drawn from our simulations here.

 We focus on \eqref{dowellnon*}--\eqref{dowellnon3*} and make clear distinctions between linear dynamics $\sigma=\iota=0$, {\bf \small [NL Stiffness]} only dynamics ($\sigma=1,~\iota=0$), and fully nonlinear dynamics---with {\bf \small [NL Stiffness]} and {\bf \small [NL Inertia]} ($\sigma=\iota=1$). 

We are interested in dynamical stability properties, as well as long-time, qualitative responses of the dynamics, to: distributed pressures (via piston theory, described in the next section), and varying initial conditions. We will measure displacements of the cantilever, and we will track things like the free end displacement curves $\big(u(L,t),w(L,t)\big)$, the arc-length of the beam as a function of time, and energies (see Section \ref{energiessec}) as a function of time.

\subsection{Dynamical Driver: Piston Theory}\label{piston}
In our simulations, we seek a simple way to test the model, affect beam stability, and ``drive" the dynamics. In line with the applications relevant to cantilever large deflections, we consider a rudimentary means for simulating the flow of gas around the cantilever. Though there are various ways to consider flow-beam coupling, the simplest is to eliminate the fluid dynamic variables altogether. This has the benefit of reducing the flow-beam system to a single non-conservative beam dynamics. Such a reduction is a dramatic simplification of complex, multi-physics phenomena, but, focusing on a simple, un-coupled model allows us to perform a thorough numerical study that can be exposited straight-forwardly. (More sophisticated related, flow-structure models are certainly explored in the rigorous mathematical literature---see, e.g., \cite{survey2,springer}.)

We consider beam dynamics interacting with a potential flow. For certain flow conditions, the dynamic pressure on the surface of the beam, $p(x,t)$, can be approximated point-wise in $x$ by an expression written in the {\em down-wash} of the fluid $W=(\partial_t+U\partial_x)w$, where $w(x,t)$ is the transverse displacement of the beam, and $U$ is the unperturbed axial flow velocity. This results in a nonlinear expression \cite{pist2} in $W$ that is linearized to produce the {\em piston-theoretic} pressure $p(x,t)$ on the beam \cite{vedeneev}:
\begin{equation}\label{linpist}
p(x,t)=p_0(x)-\beta[w_t+Uw_x].
\end{equation}
Above, $p_0(x)$ is a static pressure on the surface of the beam (for the numerical portion of this paper we will take $p_0(x)=0$). The parameter $\beta>0$ is a fluid density parameter.\footnote{See \cite{vedeneev} for a discussion of the flow non-dimensionalization, and further discussion of characteristic parameter values.}
 We consider both positive and negative values for $U$, corresponding to axial flow from clamped to free end ($U>0$---flag-like configuration \cite{huang,flag}), as well as from free to clamped end ($U<0$---inverted flag configuration \cite{inverted1,amjad}). Note that the presence of aerodynamics provides both a stabilizing term---weak damping---scaled by $\beta>0$, as well a destabilizing non-conservative term scaled by $\beta U$. (See \cite{HHWW} for more discussion.)
 
With \eqref{linpist} providing our dynamic driver, we can consider a simple non-conservative dynamics that can give rise to instability in our model, resulting in large deflections that ``test" the inextensible nonlinear effects.

	\subsection{Modal Dynamics}\label{modal}
Modal analysis, here, refers to a Galerkin method, based on the above weak formulation \eqref{weakform}, whereby solutions are approximated by in vacuo structural eigenfunctions (modes) (e.g., \cite{modalcant,vedeneev,inext1}). Since the eigenfunctions of standard elasticity operators form a basis for the state space,  a good well-posedness result for the full system justifies this type of approximation. 
This type of approximation can be dynamic, as in reducing an evolutionary PDE to a finite dimensional system of ODEs by truncation, or it can be stationary, reducing the problem of dynamic instability (for linear dynamics) to an algebraic equation. 

\subsubsection{Cantilever Modes}
Critical to any modal analysis---see, for instance, \cite{modalcant}---are the in vacuo modes (eigenfunctions) associated to the configuration. We are working with the linear Euler-Bernoulli cantilever as our approximants in $H^2_*$, and the modes and associated eigenvalues can be computed in an elementary way. These functions are {\em complete} and {\em orthonormal} in $L^2(0,L)$, as well as complete and orthogonal in $H_*^2(0,L)$ (with respect to $(\cdot,\cdot)_{H^2_*}$).

The cantilever mode shapes of interest are:
$$s_n(x)\equiv c_n(\cos(\kappa_nx)-\cosh(\kappa_nx))+ C_n(\sin(\kappa_nx)-\sinh(\kappa_nx)),$$
where the $C_n$ are obtained  by solving the associated characteristic equation: ~$\cos(\kappa_nL)\cosh(\kappa_nL)=-1.$
We have
\[ C_n=\dfrac{-c_n\big(\cos(\kappa_nL)+\cosh(\kappa_nL)\big)}{\sin(\kappa_nL)+\sinh(\kappa_nL)},\]
and the $c_n$ values are chosen to normalize the  functions in the $L^2(0,L)$ sense.

The mode numbers $\kappa_nL$ are obtained by numerically solving the characteristic equation.  

\begin{table}[b]
\begin{center}
\begin{tabular}{|l|c|}\hline
$n$& $k_nL$ \\\hline
1 & 1.8751 \\
2  & 4.6941\\
3  & 7.8548 \\
4  & 10.9955\\
5   & 14.1372\\
6  & 17.2788\\\hline
\end{tabular}
\end{center}
\caption{First 6 mode numbers for the cantilever (Clamped-Free, {\bf CF}) configuration.}
\label{table1}
\end{table}

\begin{center}
\includegraphics[width=\linewidth]{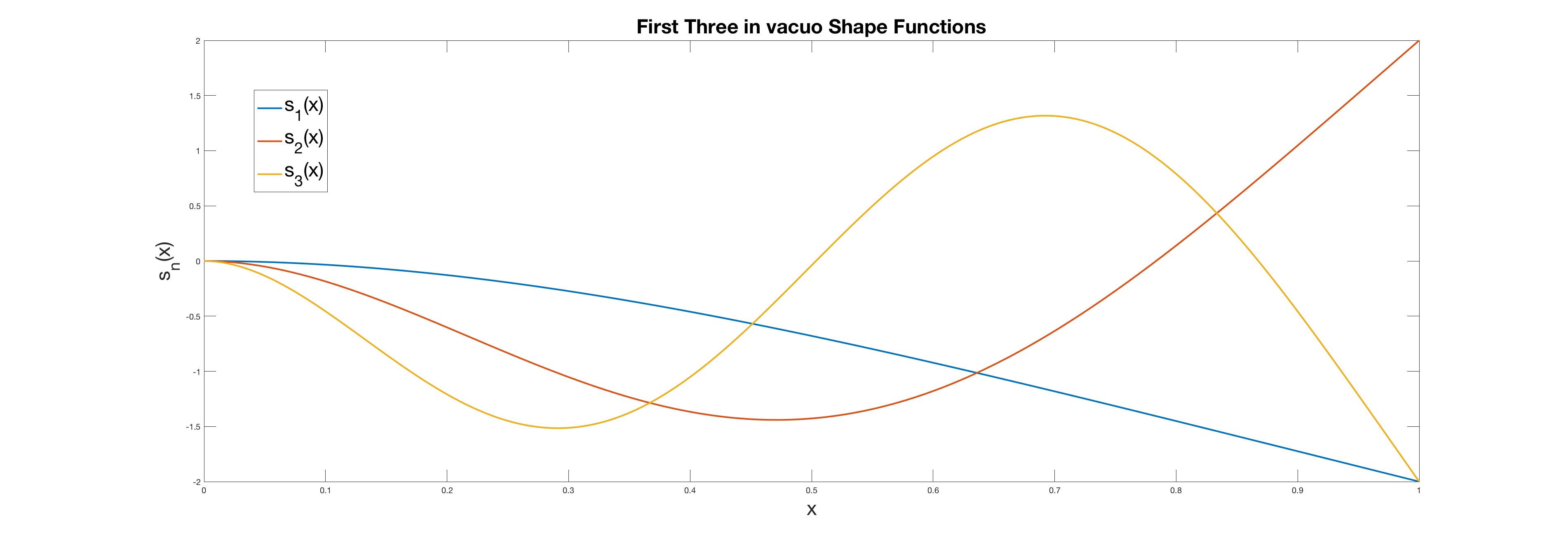}
\end{center}

\subsubsection{Calculating the Flutter Point: Reduction to Perturbed Eigenvalue Problem}\label{modalapproach}

Let us consider the Galerkin procedure for the full linear  beam equation with linear piston theory and the possibility of imposed weak damping $k_0\ge 0$: \begin{equation}\label{fornumerics}w_{tt}+Dw_{xxxx}+(k_0+\beta)w_t=-\beta Uw_x\end{equation} on $(0,L)$, with cantilever boundary conditions. We expand the solution  via the in vacuo  mode functions $\{s_j\}$ as
$w(t,x) = \sum q_j(t)s_j(x)$.  The $\{q_j\}$ represent smooth, {\em time-dependent coefficients}. Plugging this representation into the \eqref{fornumerics}, multiplying by $s_n$, and integrating over $(0,L)$ for each $n$ we obtain:
\begin{align}
\sum_m\Big[\big[q''_{m}(t)+(\beta+k_0)q'_{m}(t) +Dk_m^4q_m(t)\big](s_m,s_n)+\beta U(\partial_xs_m,s_n)q_m(t)\Big]=&~0,
\end{align}
with $'$ indicating $\partial_t$.

Orthonormality of the eigenfunctions can be invoked to produce {\em diagonal} terms, whereas the terms scaled by $\beta U(\partial_xs_m,s_n)$ are {\em off-diagonal} and give rise to the instability of the ODE system. 

To simply determine the stability of the problem as a function of the given parameters, we can invoke a standard engineering ansatz \cite{vedeneev,dowell} (and references therein): {\em assume simple harmonic motion} according to some dominant (perturbed) frequency $\widetilde\omega$; we allow possible contribution from all functions $s_n$ for $n=1,2,...,N$ via coefficients labeled $\alpha_n$:
\begin{equation}\label{preplug}w(t,x) \approx e^{-i\widetilde{\omega}t}\sum_{j=1}^{N}\alpha_j s_j(x),\end{equation} 
where $N$ is a chosen dimensional {\em truncation}. 
We multiply the modal equation by $s_n$, and then integration produces an eigenvalue problem in the perturbed frequency $\widetilde{\omega}$. 
With the off-diagonal entries $ (\partial_x s_m,s_n)$  in hand (for $1 \le m,n \le N$ with $m \neq n$), we compute diagonal terms $$\Omega_j(\widetilde \omega) = -\widetilde \omega^2-i(\beta+k_0)\widetilde\omega +Dk_j^4,~~j=1,2,...,N,$$ and we create the matrix for $1 \le n,m \le N$: 
\begin{equation}
A =A(\widetilde \omega)=[a_{mn}],~~\text{ with}~~a_{mn} = \begin{cases} \Omega_m &~\text{ for }~m =n \\
\beta U (\partial_xs_n,s_m) &~\text{ for }~m \neq n.
\end{cases}\label{modalmatrix}
\end{equation}

For chosen parameter values of $D, k_0, \beta, U, L$, we enforce the zero determinant condition for non-trivial solutions in $\widetilde{\omega}$:
$$\text{det}\big(A(\widetilde \omega)\big) =0,$$ and solve for $\widetilde \omega = [\widetilde \omega_{1},...,\widetilde \omega_N]^T$. The associated complex roots allows us to track the stability response of the natural modes to the {\em perturbation} terms.  This method is explored in-depth in \cite{HHWW} for beams across multiple configurations.
In \cite{HHWW}, however, this method is shown to be an accurate predictor of the onset of instability due to non-conservative piston-theoretic terms. 

In the simulations below, the modal method described above allows us to determine that for $D=L=k_0=\beta=1$, the onset of instability corresponding to the linear cantilever is $U_{\text{crit}}\approx 135$; this figure will arise repeatedly in our simulations below. We note that for $U>U_{\text{crit}}$, the linear dynamics have destabilized eigenvalue(s), and the linear dynamics (with no nonlinear elastic restoring force) will accordingly grow exponentially in time. We refer to this as the {\em onset} of instability due to the flow $U$.
 
\subsubsection{Nonlinear Modal Simulations}

Now, as in Section \ref{modalapproach}, let us expand the solution to the nonlinear problem \eqref{dowellnon*} as $\ds w=\sum_i s_iq_i$, where again $s_i(x)$ are the in vacuo cantilever mode shapes, with and $q_i(t)$ being smooth functions of time. Plugging the solution into the weak form \eqref{weakform} gives us a corresponding ``matrix" system in $\{q_t(t)\}$ by subsequently testing with $\phi=s_j$.

We define the following four-tensors (corresponding respectively to {\bf \small [NL Stiffness]} and {\bf \small [NL Inertia]}):
\begin{align}
\mathcal S_{ijkl}=&~(\phi_{i,xx}\phi_{j,xx},\phi_{k,x}\phi_{l,x})\\
\mathcal I_{ijkl}=&~\left(\int_0^x\phi_{i,x}\phi_{j,x},\int_0^x\phi_{k,x}\phi_{l,x}\right).
\end{align}

\begin{remark}
The following calculation for the inertial tensor connects $\mathcal I_{ijkl}$ back to the weak form \eqref{weakform}:
\begin{align*}
\mathcal I_{ijkl} = &~\left(\int_0^x\phi_{i,x}\phi_{j,x},\int_0^x\phi_{k,x}\phi_{l,x}\right) \\
 =& ~-\int_0^L \left[\left(\partial_x\int_x^L \int_0^{\xi}\phi_{i,x}\phi_{j,x}d\xi_2d\xi\right)\int_0^x\phi_{k,x}\phi_{l,x}d\xi\right]dx \\
  =&~\int_0^L \left[\left(\int_x^L \int_0^{\xi}\phi_{i,x}\phi_{j,x}d\xi_2d\xi\right)\phi_{k,x}\phi_{l,x}\right]dx.
\end{align*}
\end{remark}

Employing Einstein notation, so that $q_is_i$ is interpreted as the sum, we have the following separated form of \eqref{dowellnon*}--\eqref{dowellnon3*} taken with the piston-theoretic RHS \eqref{linpist}: {\small
\begin{align}
q''_{i}(s_i,s_j)+\left[q_i''(q_i)^2+(q_i')^2q_i\right]\mathcal I_{iiij}+\beta q_i'(s_i,s_j)+Dq_i\left[k_i^4(s_i,s_j)\right]+Dq_i^3\left[\mathcal S_{iiij}+\mathcal S_{jiii}\right]=\beta U(\partial_xs_i,s_j).
\label{odesys}
\end{align}}
This form, \eqref{odesys}, constitutes the bi-infinite {\em modal form} of \eqref{weakform} with $\iota=\sigma=1$ and $k_2=0$.
\begin{remark}
Note the {\em temporally} quasilinear term $q_i''(q_i)^2\mathcal I_{iiij}$, which may slow down time-stepping computations, as the equations are algebraically implicit in $q''$. The spatial nonlinearity---cubic type, and quasilinear---can be seen in the terms involving the tensor $\mathcal S$.
\end{remark}

The summation in equation \eqref{odesys} is truncated to include just $N$ mode functions after which we conduct a reduction of order. The resulting $2N \times 2N$ system of ODEs is solved using the \texttt{ode15i} function in MATLAB, which requires the ODE to be in the form $f(y,y_t) = 0$. To expand the summations, we use a Mathematica script and then a Python script to convert the output to valid MATLAB syntax. For the computation of the stiffness tensor components $\mathcal{S}_{ijkl}$, the inbuilt function $\texttt{integral}$ is used. 
The inertial tensor components of $\mathcal{I}_{ijkl}$, ~$ \int_0^x s_{i,x} s_{j,x},~~~\int_0^x s_{k,x} s_{l,x},$ are computed using the inbuilt $\texttt{integral}$ function and the final integral for the inner product is taken using Simpson's Rule. Once the ODE system is numerically solved, the final solution is computed by taking the corresponding linear combinations of the mode functions. Visual/graphical output can be produced as $(X,Y)=(u(x,t),w(x,t))$, where $u$ is obtained from $w$ through the effective inextensibility relation \eqref{dowellnon3*}.

\subsection{Qualitative Analysis of Numerical Simulations}\label{numerics}
For the simulations presented in this section we take the following conventions:
\begin{itemize}
\item The flags $\iota,\sigma$ take values of 0 or 1, depending on what is being discussed.
\item Non-central parameters are taken to unity: $L=\beta=D=1$.
\item The stationary pressure, $p_0(x)$ in \eqref{linpist}, is taken identically zero.
\item Imposed damping is taken to be zero, i.e., $k_0=k_2=0$.
\item Unless stated otherwise, the number of modes used in each simulation was $N=6$.
\end{itemize}

For these conventions, we mention that \cite{HHWW} and \cite{HTW} discuss piston-theoretic and structural parameter values in more depth. It is worth commenting that, even when $\iota=1$, {\em we do not invoke any damping} in our simulations below. In our theoretical results, we recall that $k_2>0$ is necessary for the existence proof to obtain when $\iota=1$ (Theorem \ref{withiota}). We choose to focus these preliminary numerical simulations on the undamped case to understand the essence of the nonlinear effects. In particular, these results below indicate precisely how the {\bf \small [NL Inertia]} effects produce issues. 

Lastly, we now specify our initial data repository for simulations below:
\begin{itemize}
\item \textbf{[1st Mode ID]} $w(0,x)=s_1(x) =  [\cos(\kappa_1x)-\cosh(\kappa_1x)]-\mathcal C_1[\sin(\kappa_1x)-\sinh(\kappa_1x)], ~w_t(0,x)=0$,  \vskip.05cm
where $\kappa_1\approx 1.8751$ is the first Euler-Bernoulli cantilevered mode number (with $L=1$) and $\mathcal C_1 \approx 0.7341$.
\item \textbf{[2nd Mode ID]}
$w(0,x)=s_2(x) =  [\cos(\kappa_2x)-\cosh(\kappa_2x)]-\mathcal C_2[\sin(\kappa_2x)-\sinh(\kappa_2x)], ~w_t(0,x)=0$,  \vskip.05cm
where $\kappa_2\approx 4.6941$ is the second Euler-Bernoulli cantilevered mode number (with $L=1$) and
 $\mathcal C_2 \approx 1.0185$.

\item \textbf{[Polynomial ID]}\quad
$w(0,x)=-4x^5+15x^4-20x^3+10x^2$, ~$w_t(0,x)=0$; 

\item  \textbf{[Linear IV]}\quad $w(0,x)=0$, $~w_t(0,x)=ax$, where the parameter $a>0$ will be increased in size as a mechanism to increase the ``size" of the initial data (in the sense of $L^2(0,1)$).
\end{itemize}

\subsubsection{Conservation of Arc Length in Numerical Simulations}
Since the inextensible dynamics are predicated on enforcing the inextensibility constraint, we posit that arc length should be approximately conserved throughout dynamic deflections. However, as we are enforcing an {\em effective} inextensibility constraint \eqref{dowellnon3*}, we expect that approximation of the full constraint \eqref{fullinext} produces errors that can be exaggerated by larger and larger deflections.

\begin{figure}[H]
\centerline{
\includegraphics[scale=0.35]{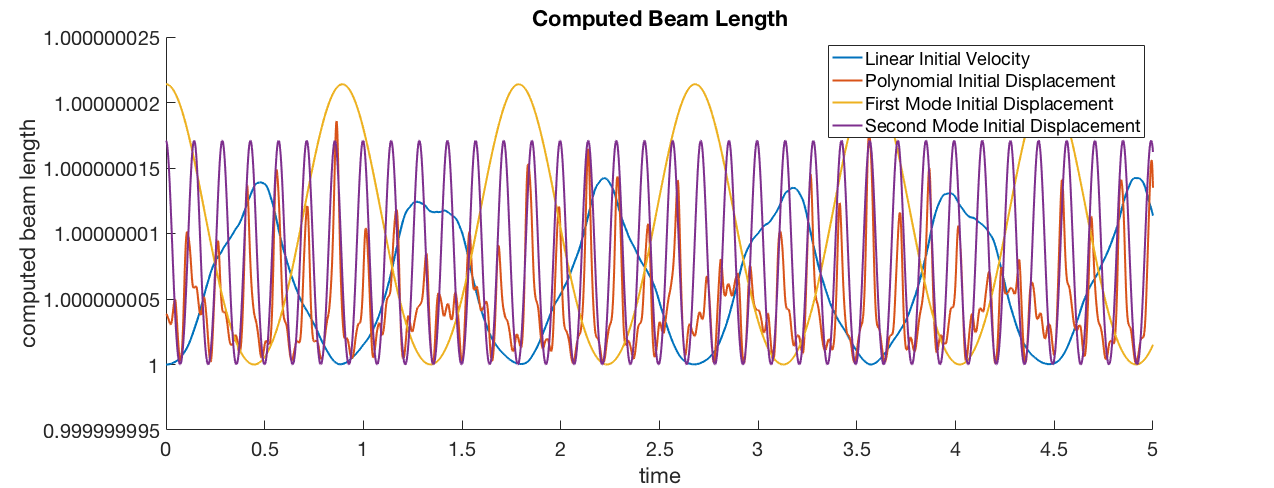}}
\caption{In vacuo computed arc length, varying initial conditions.}
\label{arclength}
\end{figure}

First, in Figure \ref{arclength}, we demonstrate that arc-length is faithfully conserved throughout deflection, across the varying initial conditions {\em for the in vacuo case} ($U=0$, $\beta=0$). These plots take active stiffness and inertia---$\iota=\sigma=1$. In these simulations, the initial velocity multiplier $a$ is set to 1.

\begin{figure}[H]
\centerline{
\includegraphics[scale=0.35]{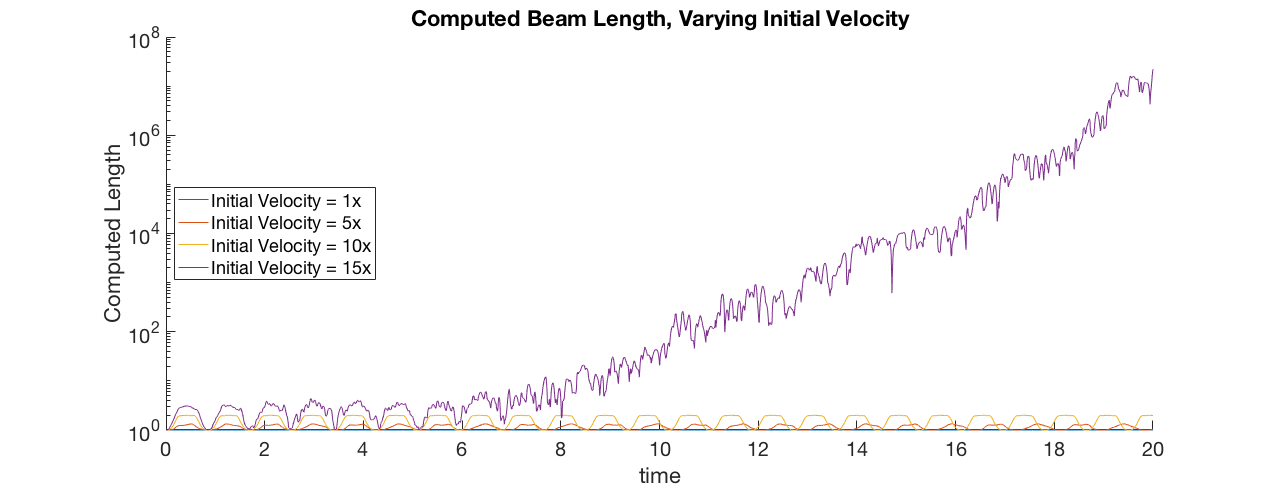}}
\caption{In vacuo computed arc length, varying initial velocity multiplier $a$.}
\label{arclengthinitvel}
\end{figure}

However, for {\bf [Linear IV]} initial conditions, increasing values of the initial velocity multiplier $a$ (with zero initial displacement) yields the degradation of arc-length conservation, as seen in  Figure \ref{arclengthinitvel}.

\begin{figure}[H]
\centerline{
\includegraphics[scale=0.35]{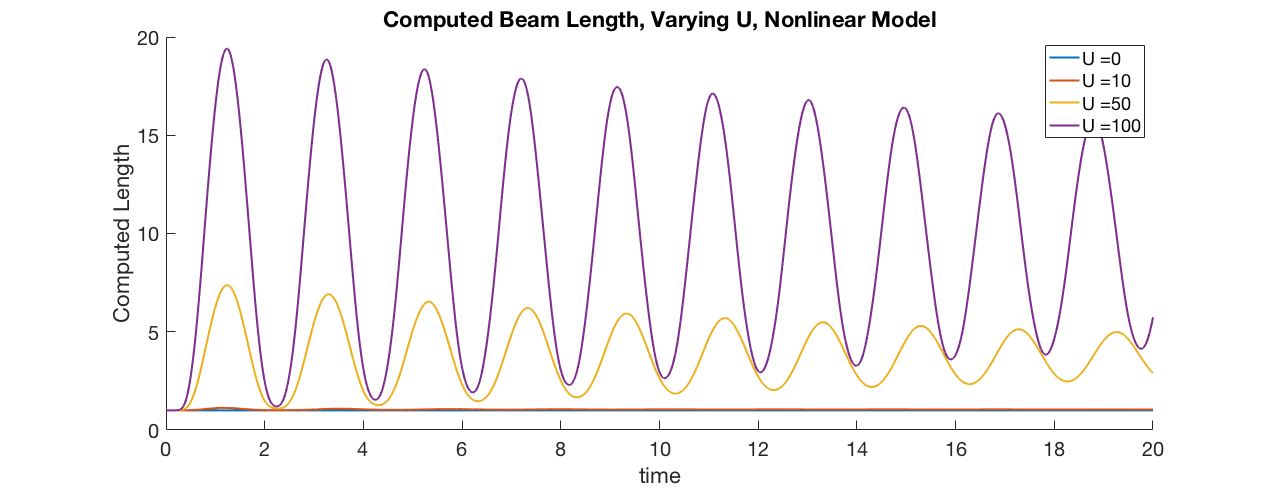}}
\caption{Full nonlinear model computed arc length, varying $U$, $\beta=1$.}
\label{arclengthU}
\end{figure}

We also see degradation in the conservation of arc length when the piston-theoretic flow is active.  Figure \ref{arclengthU} gives the computed arc length for the full nonlinear beam $\sigma=\iota=1$ for varying values of $U>0$.  The initial condition is {\bf [Linear IV]} with $a=1$.
We note that as $U$ increases, beam deflections are larger (owing to a stronger forcing), resulting in observed degradation. All flow velocities $U$ here are below the linear {\em onset} critical value of $U_{\text{crit}}\approx 135.9$. The reason for this will become clear in the discussions that follow.

\subsubsection{Computed Total Energies}
In a similar capacity to the previous section, we compute the total (nonlinear) energies associated to various situations. We are interested in tracking the temporal evolution of $E(t)$: when in vacuo dynamics are considered ($U=0$), we expect conservation of energy. When $U \neq 0$, we expect that energies will evolve, owing to the non-conservative flow effects of \eqref{linpist}.

\begin{figure}[H]
\centerline{
\includegraphics[scale=0.35]{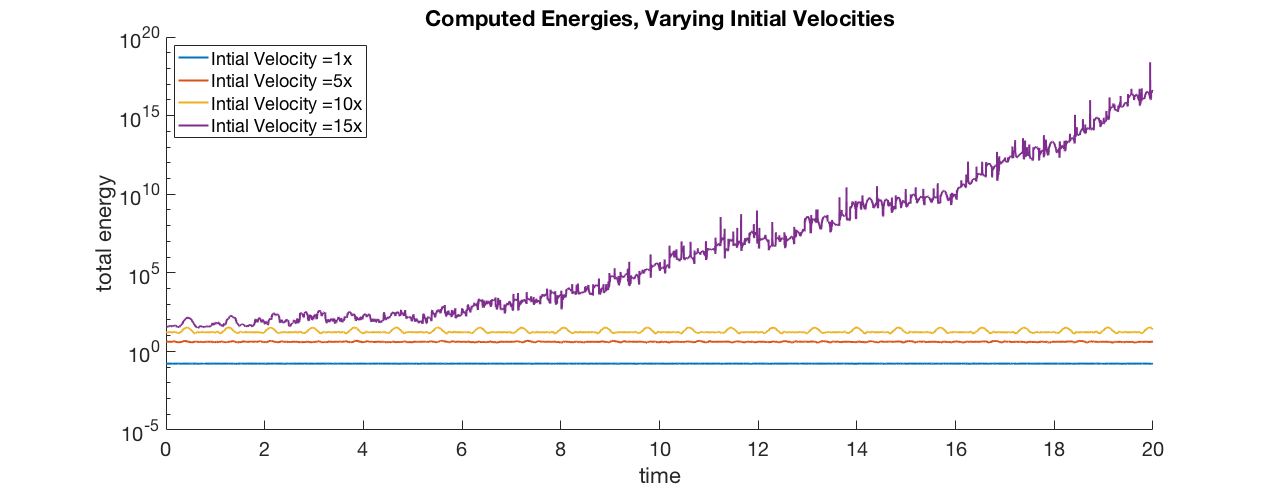}}
\caption{In vacuo total energies with {\bf [Linear IV]} varying $a$.}
\label{energyinitvel}
\end{figure}

We first examine the computed total energies for the in vacuo, fully nonlinear beam ($\iota=\sigma=1$), with varying initial velocity size in Figure \ref{energyinitvel}.
When the size of the initial data is sufficiently small, we see near perfect conservation of energy, $E(t)=E(0)$, perhaps with slight periodic effects. However, when the initial data size is large, we see that energy conservation is lost.

\begin{figure}[H]
\centerline{
\includegraphics[scale=0.35]{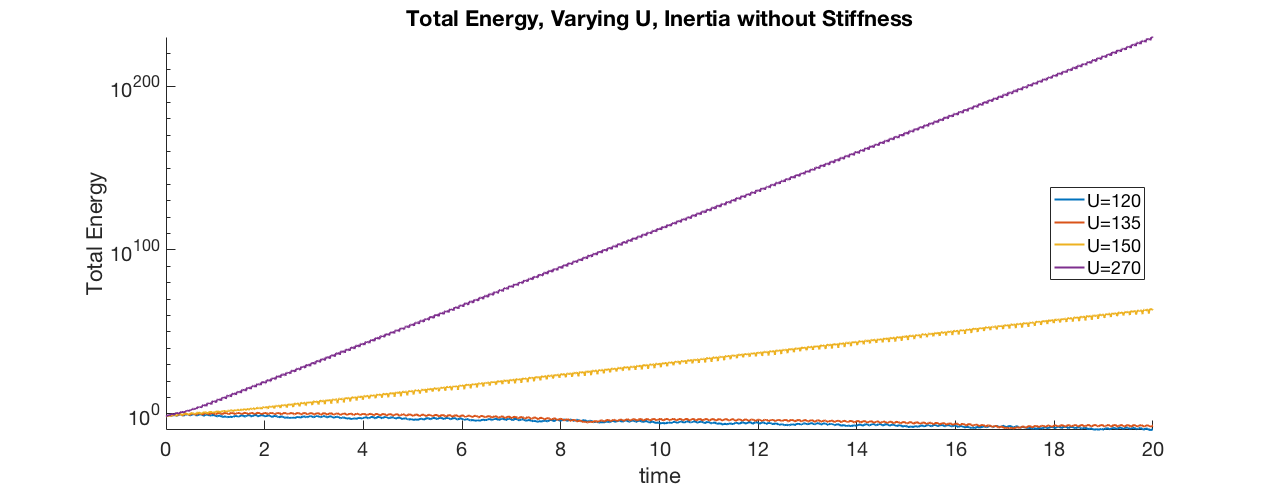}}
\caption{Total energies for the beam, varying $U$, $\beta=1$, linear model.}
\label{linearmodel}
\end{figure}

For the linear model, with stiffness and inertia turned off ($\sigma=\iota=0$), the in vacuo critical value is $U_{\text{crit}}=135.9$. Figure \ref{linearmodel} gives the computed energies for varying $U$ values below and above $U_{crit}$.
We note that for $U>U_{\text{crit}}$, we observe exponential growth in time of energies, as expected \cite{HHWW}. Below $U_{\text{crit}}$ we see exponential decay, as a function of the presence of damping in linear piston theory \eqref{linpist}. We contrast this picture with the same simulations, with active nonlinear restoring forces. (See \cite{HHWW} and \cite{HTW} for more in-depth study and discussion when an extensible beam is being considered.)

\begin{figure}[H]
\centerline{
\includegraphics[scale=0.35]{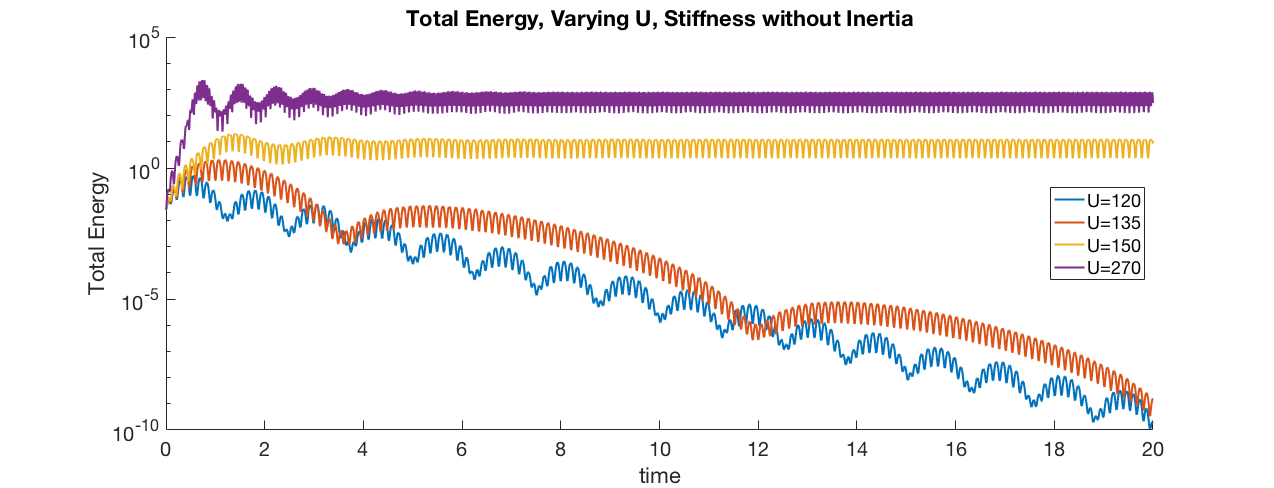}}
\caption{Energies for the beam, varying $U$,  with $\beta=1$, $\sigma=1$, and $\iota=0$.}
\label{stiffonly}
\end{figure}

First, we include stiffness only $(\sigma=1,~\iota=0)$. The energy is modified accordingly---Section \eqref{energiessec}.  Figure \ref{stiffonly} shows that the nonlinear stiffness effect is enough to provide stability in the sense that for $U>U_{\text{crit}}$, the energy plateaus. Indeed, these post-onset dynamics converge to {\em limit cycle oscillations}. Below $U_{\text{crit}}$, the trajectories which were stable in Figure \ref{linearmodel} remain so here.

\begin{figure}[H]
\centerline{
\includegraphics[scale=0.35]{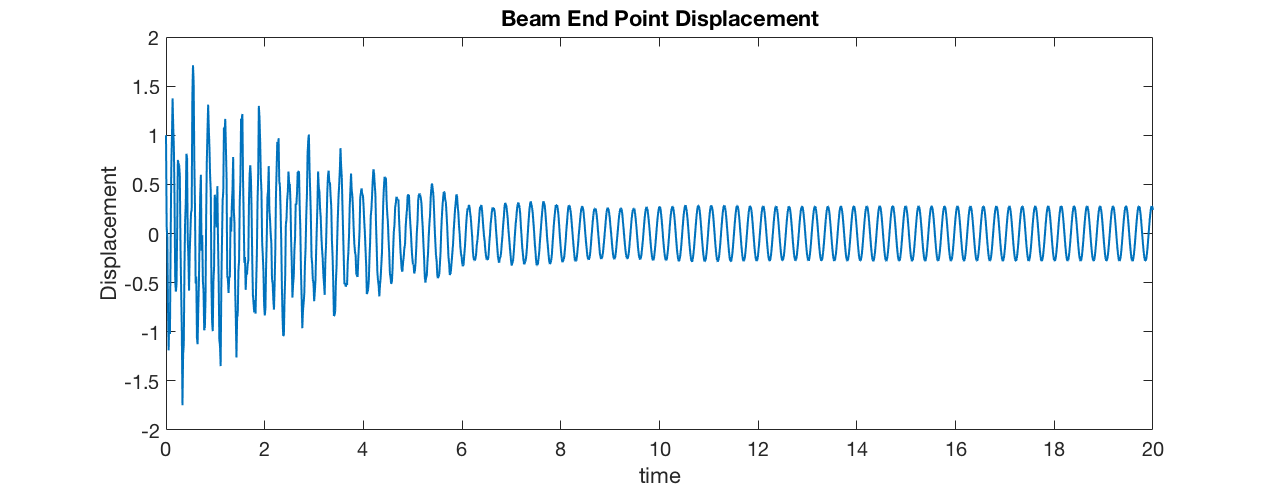}}
\caption{Limit Cycle Oscillation of beam vertical displacement, stiffness only, $U=140$, $\beta=1$.}
\label{lco}
\end{figure}

We demonstrate a limit cycle oscillation in the post-onset regime $U>U_{\text{crit}}$ for stiffness only dynamics ($\sigma=1,~\iota=0$) when $U=140$.  Figure \ref{lco} shows the beam vertical end point displacement for $U=140$.

The next question which naturally arises is: What happens to these stability (energy or displacement) plots when {\bf \small [NL Inertia]} is present in the model? More specifically, we can ask two questions: (i) Does the presence of {\bf \small [NL Stiffness]} and/or {\bf \small [NL Inertia]} affect the linear critical onset value $U_{\text{crit}}$? (ii) In the post-onset regime, what do fully nonlinear ($\sigma=\iota=1$) dynamics look like? For (i), we defer this complex question to future work, but we note that it does appear that the presence of {\bf \small [NL Inertia]} does affect---lower---the critical value for instability, however this affect is highly dependent upon initial configuration; an interesting, if not wholly surprising, observation.

\begin{figure}[H]
\centerline{
\includegraphics[scale=0.35]{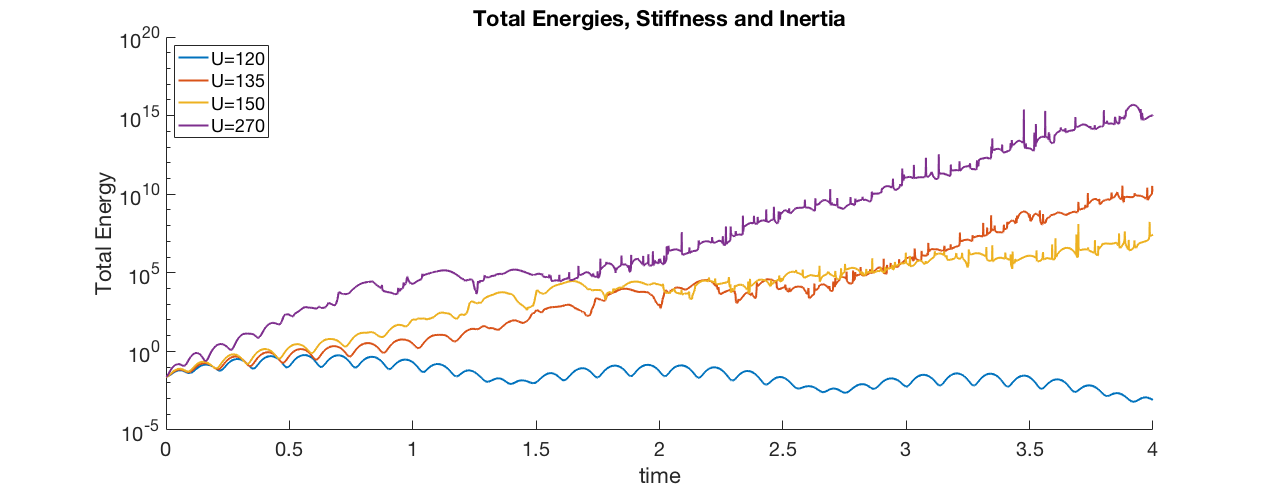}}
\caption{Total energies for the beam, varying $U$, $\beta=1$, stiffness with inertia.}
\label{stiffinertia}
\end{figure}

Indeed, with $\iota=1$, no consistent limit cycle oscillation behavior could be observed through the linear piston-theoretic RHS. This is consistent with engineering literature \cite{McHughIFASD2019,kev}. Again, outcomes are highly dependent on initial configuration.
In Figure \ref{stiffinertia}, both the stiffness and inertia terms are included in the model ($\sigma=\iota=1$).  Note that the {\bf \small [NL Inertia]} term clearly destabilizes the computation after a sufficient amount of time, and as a result, we observe blowup of total energy for a range of $U$ values, even those well below the linear critical onset value $U_{\text{crit}}$.

\subsubsection{Inverted Flag Simulations}
Interesting questions arise when the flow direction is inverted $U<0$, yielding the so called inverted flag configuration \cite{inverted1}. 
With  flow from free to clamped end, we can ask about critical values of $U$, as well as possible end behaviors.  

In principle, the structural model is robust enough to support limit cycle oscillations in this configuration, but we were unable to observe this. Piston theory, then, seems too crude a flow approximation to capture the sophisticated aerodynamic effects involved in limit cycle oscillations for the inverted flag. Convergence to non-trivial steady states, however, was observed---sometimes referred to as buckling.

\begin{figure}[H]
\centerline{
\includegraphics[scale=0.35]{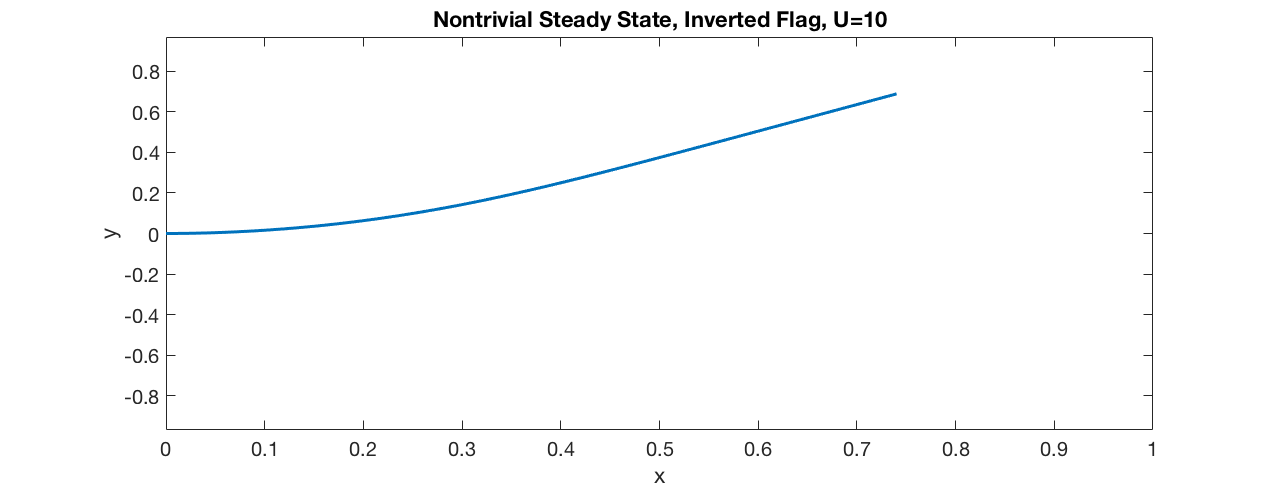}}
\caption{Nontrivial steady state displacement for the inverted flag, $U=-10$, $\beta=1$.}
\label{steady}
\end{figure}

\begin{figure}[H]
\centerline{
\includegraphics[scale=0.35]{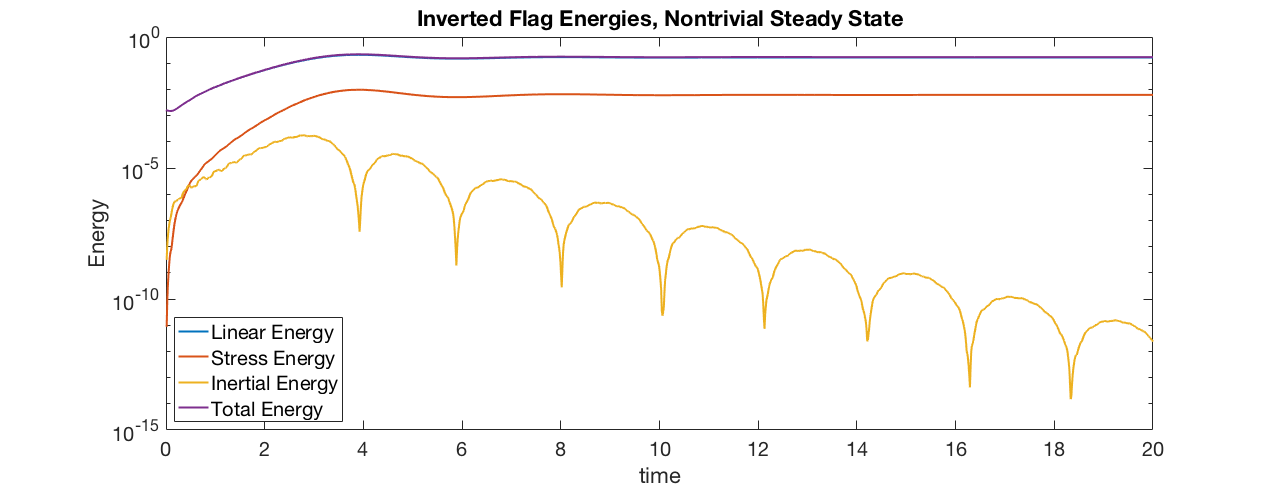}}
\caption{Nontrivial steady state energies for the inverted flag, $U=-10$, $\beta=1$.}
\label{invertedsteady}
\end{figure}

As $-U$ increases, the final steady state of the system transitions from equilibrium to a nontrivial deflected state, occurring around $\widehat{U}_{\text{crit}}=-6.3$.  Figure \ref{steady} shows a nontrivial steady state for the inverted flag configuration when $U=-10$.
Correspondingly, Figure \ref{invertedsteady} shows the different energy contributions for the inverted flag when $U=-10$.  Note the steady decay in energy associated to {\bf \small [NL Inertia]}.  As we observe convergence to a steady state, it is clear that the energy $E(t) \to const.$

\begin{figure}[H]
\centerline{
\includegraphics[scale=0.35]{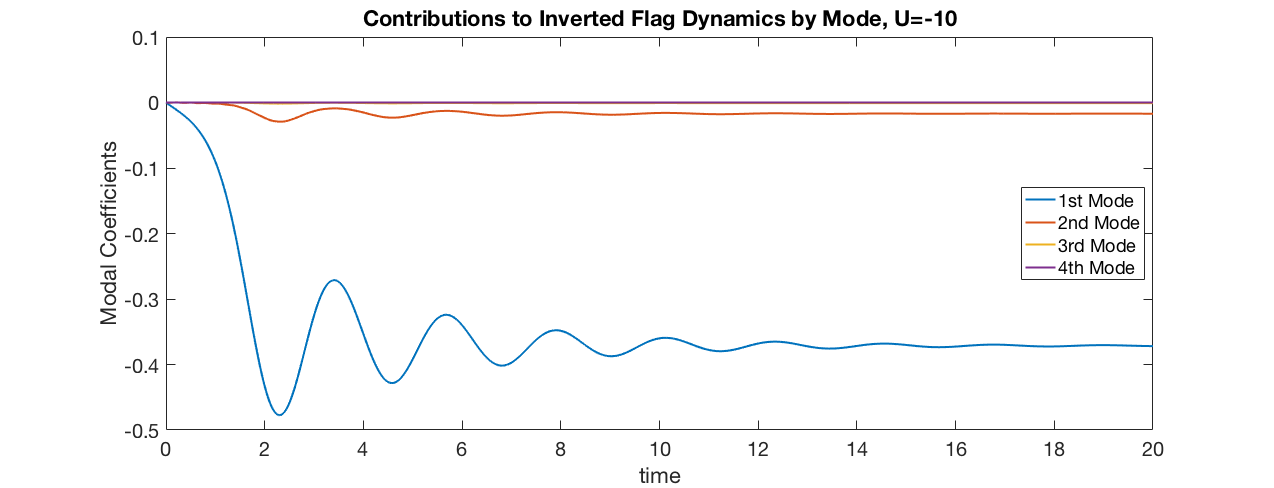}}
\caption{Modal coefficients, $U=-10$, $\beta=1$.}
\label{modalcontribution}
\end{figure}

Figure \ref{modalcontribution} tracks the values of the coefficients $q_i$ corresponding to each $s_i$ for $i=1,2,3,4$ in the inverted flag dynamics for $U=-10$.  Note that there is a contribution from both the first and second modes to the nontrivial steady state shown in Figure \ref{steady}, however there is minimal contribution from the higher modes with decay evident from the aerodynamic damping coming from \eqref{linpist}.

\begin{figure}[H]
\centerline{
\includegraphics[scale=0.35]{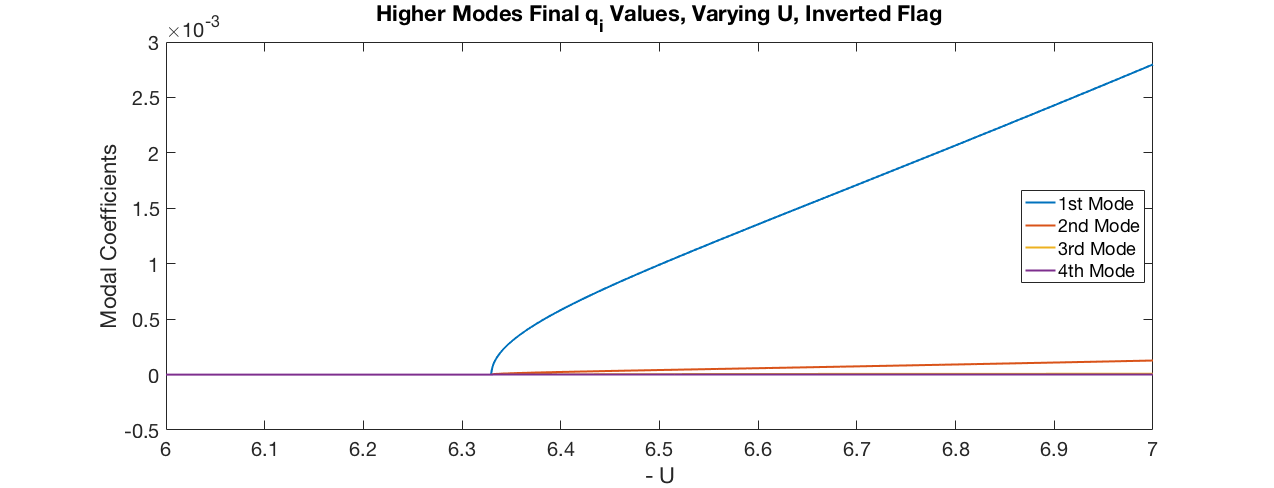}}
\caption{Modal coefficients, Varying $-6\ge U\ge -7$, $\beta=1$.}
\label{qivsU}
\end{figure}

Figure \ref{qivsU} shows how the $U$ value, ranging from $-6$ to $-7$, influences the terminal modal coefficients at $T=20$.  Around $U=-6.3$ we see the deflected steady state emerge.

\subsubsection{Influence of Number of Modes on Simulations}

Below is a plot of the total energy for $U=U_{crit}=130$, fully nonlinear model ($\iota=\sigma=1$) {\em as we increase the number of modes in the simulation}.  From the plot it is clear that varying the number of modes does effect the onset of instabilities, as the energy remains bounded for only 3 modes, but we see exponential growth with 6 modes.

\begin{figure}[H]
\centerline{
\includegraphics[scale=0.35]{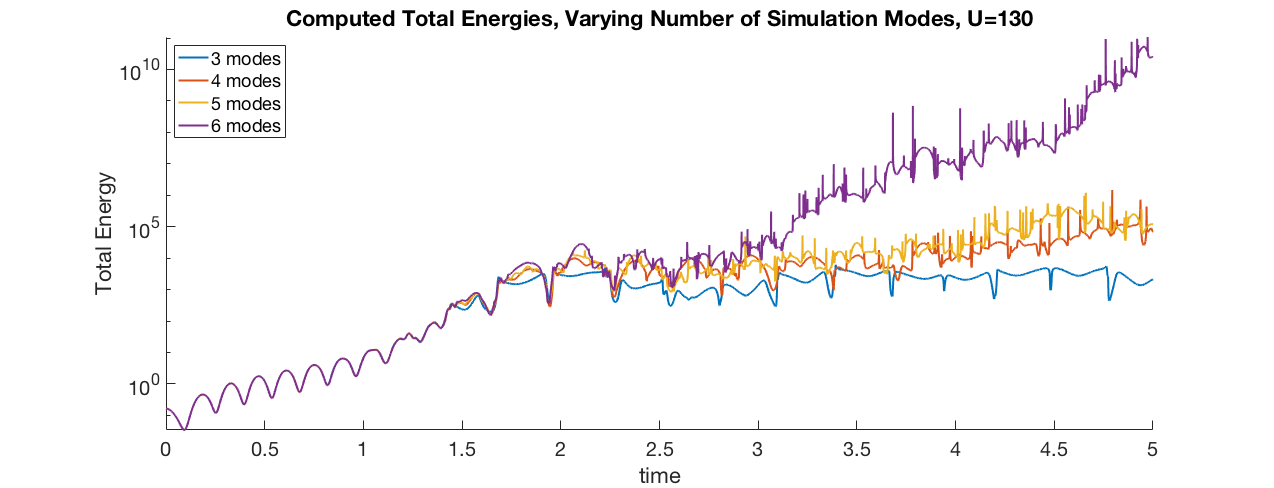}}
\caption{Energies for the fully nonlinear model, $U=135.9$, $\beta=1$.}
\label{energiesbymode}
\end{figure}

\begin{figure}[H]
\centerline{
\includegraphics[scale=0.35]{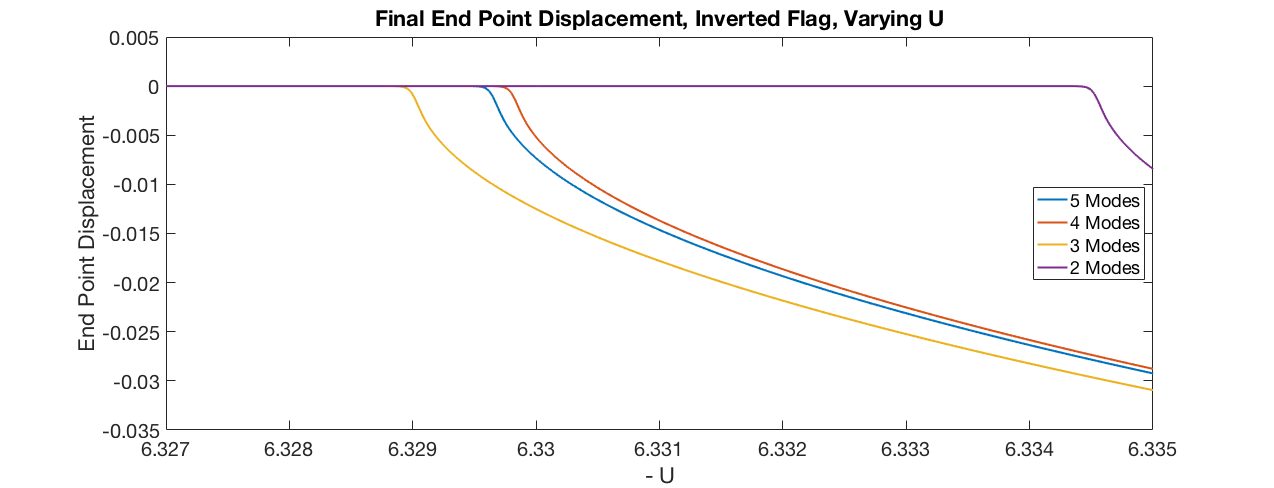}}
\caption{Final $x=L$ displacement, inverted flag, varying $U$ and number of modes, $\beta=1$.}
\label{finalbymode}
\end{figure}

For the inverted flag system ($U<0$), the number of modes  influences where the beam begins to transition to a deflected steady state, as shown in Figure \ref{finalbymode}.

\subsection{Numerical Conclusions}\label{numconc}
We briefly provide conclusions drawn from the previous section.
\begin{itemize}
\item Arc length and energy conservation are reasonably satisfied for the inextensible beam dynamics, at least for solutions with small enough deflections. Both breakdown when deflections are sufficiently large, and we expect this is due to the violation of the assumption $u_x<<1$.
\item {\small \bf [NL Stiffness]} provides a strong enough (quasilinear) restoring force to bounded post-critical trajectories and provide limit cycles. This is similar to the situation for extensible beams, where the nonlinearity is semilinear (as in Krieger or von Karman type); see \cite{HHWW,springer}.
\item {\small \bf [NL Inertia]} is the challenging term. With inertia in force,  large deflections become problematic. Specifically, simulations break or become unphysical (for instance, with the beam bending back on itself or kinking.) Considering piston-theoretic pressures, no consistent prediction of onset of instability can be given when $\iota=1$. Moreover, no consistent post-critical behavior was found; specifically, no stable limit cycles were observed when $\iota=1$. 
\item For the full nonlinear model of the inextensible beam, we see many effects in the qualitative behavior that depend on the initial data size and type. Noting that our results in Section \ref{wellpresults} are {\em local} results, this is not surprising. Moreover, the quasilinear (in time and space) nature of the model seems to suggest that dependence on data is unavoidable for theoretical and numerical results. 
\item For the inverted flag configuration $U<0$, non-trivial steady states were  obtained in a post-critical regime. Again, no limit cycles were observed when $\sigma=1$, regardless of the value of $\iota$.
\item The number of modes used in simulating the dynamics affects stability and qualitative properties. This is atypical, by engineering standards, where for structures that are non cantilevers, so called {\em modal convergence}, is observed fairly uniformly for $N\approx 4$. Owing the highly nonlinear nature of this model, as well as the large deflection nature and free boundary condition, modal convergence is not observed for such low mode numbers, and eventual behavior is dependent on $N$.
\end{itemize}

\section{Open Questions and Future Work}\label{future}
Future work, extending from the results presented here and in \cite{maria}, includes a variety of challenging modeling, analysis, and numerical problems.

\subsection{Optimal Damping}\label{optimaldamp}

As alluded to in Sections \ref{damping} and \ref{w/inertia}, many types of abstract interior damping are possible for the beam dynamics. Specifically, here, we look at damping of the form $\mathcal A^{\theta}w_t$, as described in Section \ref{damping}. For general $\theta \in (0,1)$ this is a nonlocal operator that can depend on the boundary conditions. For $\theta=0$, we have weak damping ($k_0>0$ in earlier sections) and for $\theta=1$ we have the strong damping of the form utilized here ($k_2>0$). Independent of the physical/engineering considerations discussed earlier about the physical interpretation of, for instance, $\mathcal A^{1/2}w_t$, one can ask the following mathematical question: \begin{quote} What strength of damping (what power $\theta \in (0,1)$) is sufficient to obtain energy estimates resulting in well-posedness for the inextensible beam that allows {\bf \small [NL Inertia]}? \end{quote}
Obviously $\theta=1$ is sufficient for our purposes here---taking $k_2>0$---but is this the optimal power? We conjecture, in fact, that $\theta=1/2$ corresponding to square root damping  is sufficient to obtain estimates. Secondly, what we must address is the non-equality between $\mathcal A^{1/2}w_t$ and $-\partial_x^2w_t$  ($k_1>0$ in \eqref{linearplate}), as per \cite{beamdamping}. It is clear that weak damping---$k_0>0$---is not sufficient. 

We also mention that in \cite{beamdamping,fab-han:01:DCDS} so-called indirect damping mechanisms are introduced that effectively operate as square-root type damping and do not have the issues with physical interpretation/boundary conditions in the cantilever configuration. Such damping is of interest here for providing control of the nonlinear inertia terms.

\subsection{Global Solutions and Stability}

In line with \cite{xiang} and other papers on quasilinear beam and plate equations, we seek global solutions when damping is present. Indeed, since damping is---rather oddly---necessary for us to obtain existence for the full system ($\sigma=\iota=1$) when inertia is present, we may ask what sort of stability for the system is gained as a byproduct. We expect that, the stiffness only dynamics ($\iota=0$) will have typical quasilinear behavior when strong damping is present ($k_2>0$). Specifically, we anticipate that the time of existence can be taken to be $T^*=\infty$ if the initial data is confined to a ball in the sense of $\mathscr H^s$. Such a proof depends on Lyapunov methods to show exponential decay for sufficiently small data, or sufficiently large damping coefficient $k_2$. When $\iota=1$ and $k_2>0$, this question is more delicate. At present, it is unclear whether a global-existence-with-small-data/exponential decay result will hold when {\bf \small [NL Inertia]} is active. 

\subsection{Uniqueness of Strong Solutions with $\iota=1$}

As we allude to above, no claims are made about uniqueness of strong solutions when {\bf \small [NL Inertia]} is active. We simply do not have the regularity needed---even with $k_2>0$---to obtain closed estimates on the difference of trajectories when $\iota=1$. It is clear that higher regularity of trajectories is necessary, precipitating the need to close estimates in higher topologies. This corresponds to highly complex polynomial-type differentiated versions of \eqref{dowellnon*}--\eqref{dowellnon3*}, and the associated energetic approaches.

\subsection{Other Non-Conservative Models}
From the modeling point of view, there are other beam configurations and models of interest which are of interest in engineering. These include the free-free beam, as well as cantilevers driven by piston-theoretic forces and non-conservative follower forces. The former represent linear or nonlinear aerodynamic forces that are distributed across the beam---see \cite{McHughIFASD2019,pist2,vedeneev,HHWW}; the latter represent forces that are purely tangential at a free end, throughout deflection (making the boundary force necessarily nonlinear, since it depends on the deflected state)---see \cite{follower} and references therein. In these situations, the modeling presented here (consistent with \cite{inext1}) is possibly altered via the handling of the Lagrange multiplier $\lambda$ that enforces inextensibility. Moreover, well-posedness in these situations becomes an open question once again, since the non-conservative, typically lower order terms, interact with both quasilinear and nonlocal effects resulting from inextensibility.

Beyond existence and uniqueness theory, the question of global stability or at least global existence when damping is present, is of course altered by the presence of these non conservative terms. It is then natural to ask how nonlinear effects owing to inextensibility respond to non-conservative and perhaps nonlinear effects such as piston theory or follower forces. As we see from the use of linear piston theory in Section \ref{numerics} as a driver of inextensible dynamics, many behaviors are possible. This leads to deep questions about the effect of such non-conservative terms on stability properties of the model, including time of existence, perhaps even when strong damping is present. In line with the long-time behavior analysis in \cite{HTW}, as well as the qualitative numerical analyses in \cite{HHWW}, future studies will address the stability and time of existence for the inextensible cantilever dynamics when damping size, initial data, and non-conservative coefficients are varied.

\subsection{Obtaining Limit Cycle Oscillations for Full Inextensible Dynamics}
In line with the work in the engineering literature making use of follower forces and higher order piston theory, we hope to produce limit cycle oscillations for the dynamics {\em when inertia is active}. The first set of numerical tests to be run involve the impact of strong damping $k_2>0$ on controlling the ``bad" behaviors of {\bf \small [NL Inertia]}. Beyond the effect of damping (perhaps to permit limit cycles), we will attempt to simulate limit cycles by involving the more sophisticated {\em nonlinear piston theory}; recent work \cite{kev} has indicated that large deflection piston theory produces such limit cycle oscillations. Future work will also tackle the problem of determining $U_{\text{crit}}$ with nonlinear effects active; specifically, how the linear $U_{\text{crit}}\approx 135.9$ would lower when $\sigma=1$ and/or $\iota=1$. Even formulating this problem is difficult, as the effect of nonlinearity here on stability seems to be highly dependent on initial conditions.

\subsection{Inextensible Cantilevered Plates}
Another topic of interest is the development of a  {\em 2-D  inextensible theory} for plates  (akin to extensible von Karman theory for plates \cite{koch, springer,lagnese}) defined on {$\Omega \subset \mathbb R^2$}. We note that the expression of the free boundary conditions, as well as the operator theoretic setup associated with free boundary conditions, are much more complex for plates. Stability problems and numerical analyses for cantilevered plates in axial flow have appeared in the engineering literature for some time \cite{ELSS,dowell4,inext2}, but, as with the beam, no mathematical theory seems to exist.

To provide some 2-D modeling insight, let {\small $\mathbf u=(u^1,u^2)$} represent in plane displacements; then, inextensibility implies that in-plane strains are zero, from which we approximate: 
$$0=\partial_{x_j}u^j+\frac{1}{2}\big([u^1_{x_j}]^2+[u^2_{x_j}]^2+[w_{x_j}]^2\big) \implies \partial_{x_1}u^1=-\frac{1}{2}[w_{x_2}]^2,~~\partial_{x_2}u^2=-\frac{1}{2}[w_{x_1}]^2,~~\nabla \cdot \mathbf u=-w_{x_1}w_{x_2}.$$

 The potential energy, in this case, is given by $$E_P=|| (1+|\nabla w|^2)^{1/2}\Delta w||^2_{L^2(\Omega)}.$$
From these identities, equations of motion can be derived via the same variational technique described in Section \ref{inextcant}.

In the case of rectangular plates, 2-D mode functions can be taken, roughly, as products of 1-D cantilever modes along with free-free modes. Yet more robust theoretical and computational approaches are called for. Numerically, the development of spectral and FEM methods seem particularly relevant. Analytically, the loss of the Sobolev embedding  ({$H^1(\Omega)\hookrightarrow C([0,L])$}) is critical for energy methods described here (in the analysis of differentiated equations), and thus presents a dimensionally-dependent challenge.

\section{Acknowledgements} 
The authors acknowledge the generous support of the National Science Foundation: M. Deliyianni and J.T. Webster's research contributions here were partially supported by NSF-DMS-1907620; J.Howell's research contributions here were partially supported by NSF-DMS-1908033. V. Gudibanda acknowledges support through Carnegie Mellon's SURF program.

\end{document}